\documentclass[11pt]{amsart}
\usepackage{amsmath,amsthm,amsfonts,amssymb,bbm, setspace,graphicx,float,subfigure,xcolor,mathrsfs,mathtools, Lutsko_Style}
\usepackage{todonotes}
\usepackage{enumitem}
\setlist[enumerate,1]{label={(\roman*)}}

\begin{document}

\author{Christopher Lutsko}
\address{University of Houston, Mathematics Department, Houston, Texas, 77004, USA}\email{clutsko@uh.edu}

\author{Tobias Weich}
\address{Universität Paderborn, Institut für Mathematik, Institut für photonische Quantensysteme (PhoQS), Warburger Str. 100,
33098 Paderborn, Deutschland} \email{weich@math.upb.de}

\author{Lasse L. Wolf}
\address{Universität Paderborn, Institut für Mathematik, Warburger Str. 100,
33098 Paderborn, Deutschland} \email{llwolf@math.upb.de}

\title[Polyhedral bounds and temperedness]{Polyhedral bounds on the joint spectrum and temperedness of locally symmetric spaces}

\keywords{Semisimple Lie groups, joint spectrum, temperedness}

\begin{abstract}Given a real semisimple connected Lie group $G$ and a discrete subgroup
$\Gamma < G$ we prove a precise connection between growth rates of the group
$\Gamma$, polyhedral bounds on the joint spectrum of the ring of invariant differential operators,
and the decay of matrix coefficients. In particular, this allows us to completely characterize temperedness of $L^2(\Gamma\backslash G)$ in terms of Quint's growth indicator function. As an application of our sharp polyhedral bounds we prove temperedness of $L^2(\Gamma\backslash G)$ for all Borel Anosov subgroups $\Gamma$ in higher rank Lie groups $G$ not locally isomorphic to $\mathfrak{sl}_3(\mathbb{K}), \mathbb{K}=\R,\C,\mathbb H,$ or $\mathfrak{e}_{6(-26)}$.
\end{abstract}

\subjclass[2020]{22E46, 58C40}

\setcounter{tocdepth}{2}  \maketitle 

\maketitle

\section{Introduction}

Consider a locally symmetric space $\Gamma\backslash G/K$, where $G$ is a real connected semisimple non-compact Lie group with finite center, $K$ is a maximal compact subgroup, and $\Gamma< G$ is a discrete subgroup. When the group $G$ has rank one, there is an important connection between:
\begin{enumerate}[label = (\roman*)]
  \item The bottom of the $L^2$-spectrum of the Laplace-Beltrami operator.
  \item The exponential growth rate of $\Gamma$ points in $G/K$ in a ball of growing Riemannian distance (given by the \emph{critical exponent} $\delta_\Gamma$, see \eqref{eq:def_critical_exp}).
  \item The decay rate of matrix coefficients of $L^2(\Gamma \backslash G)$ (i.e.~temperedness).
\end{enumerate}
For $G=\SL_2(\R)$ the connection between (i) and (ii) was achieved in the seminal work on the subject by Elstrodt \cite{MR360472,MR360473,MR360474} and Patterson \cite{pattersonlimitset} (see Subsection \ref{ss:rank one}). The relation between (i) and (iii) is a direct consequence of the explicit knowledge of all unitary irreducible $\SL_2(\R)$-representations and one deduces that $L^2(\Gamma\backslash G)$ is tempered if and only if $\delta_\Gamma\leq 1/2$. However, the theorem of Elstrodt-Patterson is equally of interest for $\delta_\Gamma>1/2$ as this ensures an eigenvalue of $\Delta$ below $1/4$, often called an exceptional eigenvalue. These eigenvalues determine the spectral gap for the Laplacian. As such, controlling for the location of exceptional eigenvalues plays a pivotal role in many important works. For example, the uniform spectral gap estimates for congruence subgroups and applications to expander graphs obtained by Gamburd \cite{Gamburd2002} and affine sieves by Bourgain, Gamburd, and Sarnak \cite{BGS} (see also the recent result of Calder\'on-Magee \cite{CalderonMagee2024}) and the uniform spectral gap estimates for random covers of Magee and Naud \cite{MN}. See also, the recent work of Anantharaman and Monk \cite{AM1, AM2} with regards to the spectral gap of a random surface.

The aim of this article is to prove a generalization of the  Elstrodt-Patterson theorem for the joint spectrum of invariant differential operators on higher rank locally symmetric spaces and to reproduce the above trichotomy in full generality.

Before stating the main theorem we need to establish some notation. Recall that $G$ admits a Cartan decomposition $G=K \exp(\overline{\mathfrak a_+}) K$. Hence, for every $g\in G$ there is a $\mu_+(g)\in \overline{\mathfrak a_+}$ such that $g\in K\exp(\mu_+(g))K$. $\mu_+(g)$ can be thought of as a higher dimensional distance $d(gK,eK)$. Generalizing the critical exponent $\delta_\Gamma$ to higher rank, Quint \cite{Qui02} introduced the notion of the growth indicator function $\psi_\Gamma\colon \mathfrak{a} \to \R\cup\{-\infty\}$:
\begin{align*}
		\psi_\G(H)\coloneqq \|H\| \inf_{H\in\mathcal C} \inf\left\{s\in\R\mid \sum_{\g\in \G,\mu_+(\g)\in\mathcal C} e^{-s\|\mu_+(\g)\|} <\infty\right \},
	\end{align*}
where the first infimum runs over all open cones $\mathcal C\subseteq\mathfrak{a}$ with $H\in \mathcal{C}$.
See Subsection \ref{ss:gif} for more details.
We will measure its size with respect to $\mu \in \mathfrak{a}^\ast$
by the modified $\mu$-critical exponent
\[
\delta_{\Gamma}'(\mu) \coloneqq\inf\{t\in\mathbb{R} \mid t\mu(H) > \psi_\Gamma(H) -\rho(H) \;\forall H\in \overline{\mathfrak{a}_+}\}
\]
which equals the abscissa of convergence for the series
\[
\sum_{\gamma\in \Gamma} e^{-(s\mu +\rho)(\mu_+(\gamma))}
\]
by \cite[Prop~3.1.8]{Qui02} 
(see also  \cite[Lemma~2.1]{WZ23} for a statement closer to the notation of the present paper).

Note that the abscissa of convergence of
\[
  \sum_{\gamma\in \Gamma} e^{-s\mu(\mu_+(\gamma))}
 \]
 is a widely used quantity in the study of discrete groups, that goes back at least to the work of Quint \cite{Qui02}  and is often called the $\mu$-critical exponent $\delta_\Gamma(\mu)$ of $\Gamma$. Hence, we have chosen the name \emph{modified critical exponent} because in our setting the $\rho$-shift naturally occurs (see e.g. \cite{OMT23, KMO24} for more recent work on these critical exponents).

In higher rank, the role of the Laplacian is played by the
full algebra of invariant differential operators on $G/K$ which we denote by $\mathbb D(G/K)$. It is convenient to parametrize the joint spectrum of this algebra via the Harish-Chandra isomorphism by a $W$-invariant subset
$\widetilde \sigma_\Gamma \subseteq \mathfrak{a}^\ast_\C \cong \C^{\text{rank}(G/K)}$ (see \cite{DKV, JL} for more details on this relation).
In general,
\begin{align}
	\label{eq:generalboundspectrum}
  \Re \widetilde \sigma_\Gamma\subseteq \operatorname{conv}(W\rho),
\end{align}
where $\rho$ denotes the usual half-sum of positive restricted roots and $\operatorname{conv}(W\rho)$ is the polyhedron described by the convex hull of the Weyl orbit of $\rho$
(see Section~\ref{sec:spherical_dual}).
Moreover, $\wt \sigma_\Gamma \subseteq  \{\lambda\in \mathfrak{a}_\C^\ast \colon -\overline{\lambda}\in W\lambda\}$
so that $\Re \wt \sigma_\Gamma \subseteq \mathfrak{a}^{\ast,\mathrm{Her}}\coloneqq
\{\lambda\in \mathfrak{a}^\ast \colon -\lambda\in W\lambda\}$.

\begin{figure}
	\centering
	\includegraphics[width = \textwidth, trim = 0 3.5cm 0 0cm, clip]{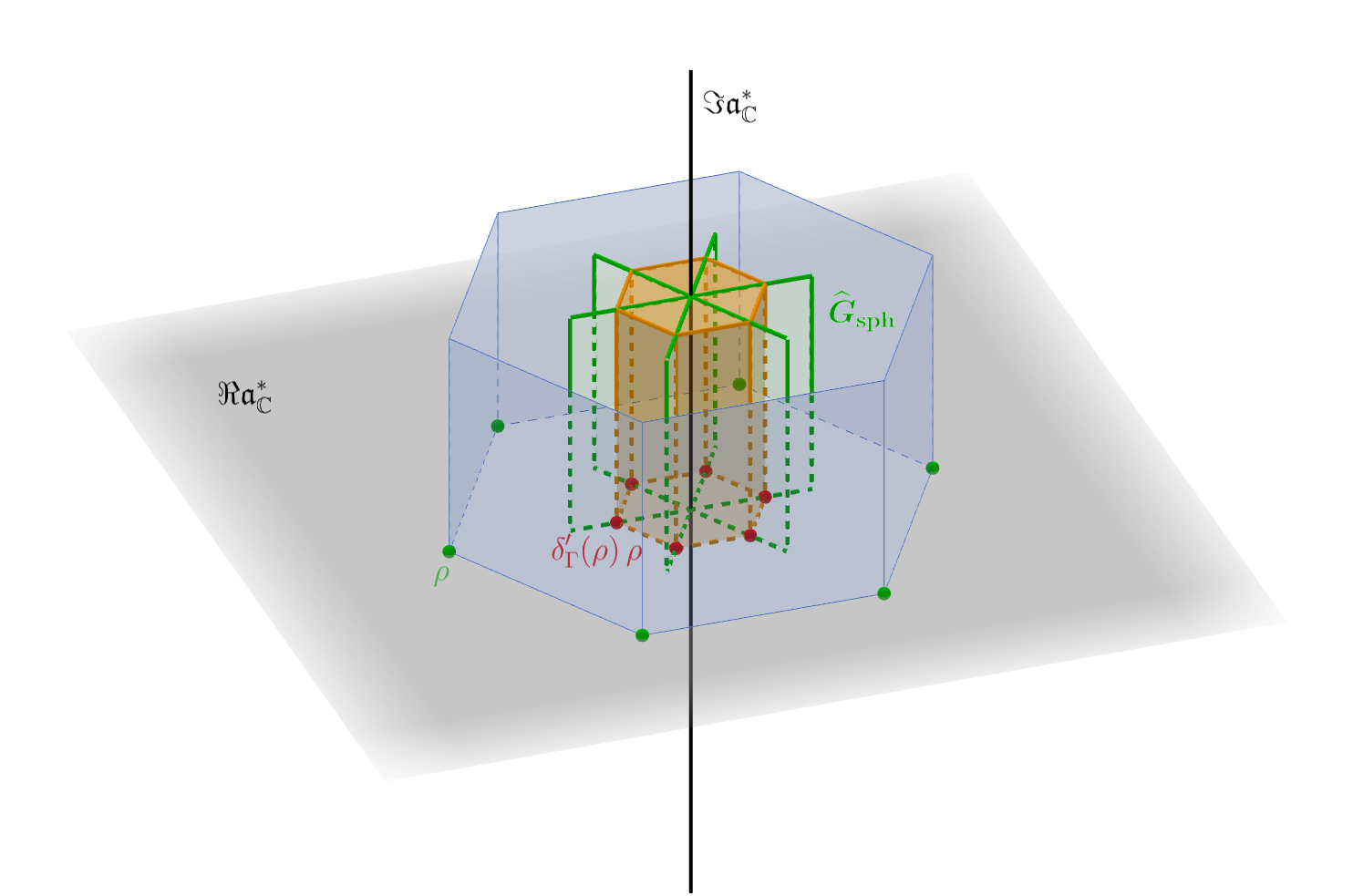}
	\caption{Visualization for $G=\mathrm{SL}_3(\R)$ in the case $\delta_\Gamma'(\rho)\geq 0$.
The gray plane is the real part of $\mathfrak{a}_{\C}^\ast$.
The two-dimensional imaginary part is depicted as a one-dimensional $z$-axis.
The green planes together with $W\rho$ is where the joint spectrum can actually occur,
i.e.~this is $\widehat G_{\mathrm{sph}}$. The blue hexagonal tube is the region $\{\Re \lambda \in \operatorname{conv}(W\rho)\}$ which is the general bound \eqref{eq:generalboundspectrum} for the real part of the joint spectrum. The orange tube is the restricted region containing $\widetilde \sigma_\Gamma$ by Theorem~\ref{thm:main}.
By Theorem~\ref{thm:main} we know that there is spectrum arbitrarily close to the boundary of the orange tube.
Proposition~\ref{prop:sl3} shows that this occurs actually at $\delta_{\Gamma}'(\rho)\rho$ (red).
We note that $\mathfrak{a}^{\ast,\mathrm{Her}}=\R\rho$ and therefore Theorem~\ref{thm:main} can only be applied to $\mu=\rho$.
}
	\label{fig:sl3}
\end{figure}

Furthermore, we introduce the \emph{polyhedral norm} which is the key ingredient to formulate our main theorem:
	For any linear functional $\lambda : \mathfrak{a} \to \R$, and any $\mu\in\mathfrak a^*$ as abve,
\[
	\|\lambda \|_{\mathrm{poly},\mu} = 
	\sup _{w\in W,H\in \overline{\mathfrak{a}_+}} \frac{w\lambda(H)}{\mu(H)}.
        \]
The terminology polyhedral norm stems from the fact that this is a vector space norm on $\mathfrak a^*$ whose balls are polyhedra spanned by the Weyl translates of $\mu$, i.e.
\begin{align*}
	\{\lambda \in \mathfrak{a}^\ast, \ \|\lambda\|_{\mathrm{poly},\mu} \le R\} = R \operatorname{conv}(W \mu).
\end{align*}
Thus, the general bound \eqref{eq:generalboundspectrum} on the joint spectrum is equivalent to saying that, for arbitrary $\Gamma$,  $\|\Re \lambda\|_{\mathrm{poly},\rho}\le 1$ for all $\lambda \in \wt{\sigma}_\Gamma$ (cf. Figure~\ref{fig:sl3} for a visualisation for $\SL_3(\R)$).
We also note that
\[
	\delta_{\Gamma}'(\mu) = \sup_{w\in W, H\in \overline{\mathfrak{a}_+^\ast}}
	\frac{\psi_\Gamma(w^{-1}H)-w\rho(H)}{\mu(H)}
\]
where we let $W$ act on $\mathfrak{a}$ by duality.
Hence, $\delta_{\Gamma}'(\mu)$ could be thought of the polyhedral norm
(with respect to $\mu$) of the positively homogeneous function $\psi_\Gamma -\rho: \mathfrak a \to \R$.
However, as it doesn't define a norm on the space of positively homogeneous functions $\mathfrak{a}\to \R$,
we choose a different notation.

As a last ingredient let us introduce the exponential decay rate of matrix coefficients:Recall that $L^2(\Gamma\backslash G)$ is a unitary representation and $L^2(\Gamma\backslash G)^K$ denotes the $K$-invariant vectors that are in 1:1 correspondence to $L^2(\Gamma\backslash G/K)$ and are the natural vectors to study in the context of the spectral theory of $\mathbb D(G/K)$ on $L^2(\Gamma\backslash G/K)$ (cf. Proposition~\ref{prop:joint_spec}).
	Let $\theta_\Gamma(\mu) \geq 0 $ denote the infimum of all $\theta'\geq 0$ such that, for all $ v\in \overline{\mathfrak{a}_+},$ and $ f_1,f_2 \in L^2(\Gamma\backslash G)^K$, and all $\varepsilon>0$, one has
\begin{align*}
	\left|\langle (\exp v)f_1,f_2\rangle_{L^2(\Gamma\backslash G)}\right|
					\leq C e^{\vep\|v\|+(\theta'\mu- \rho)(v)} \|f_1\|_2\|f_2\|_2,
			\end{align*}
		for some $C>0$ independent of the choice of $v$ or functions $f_1, f_2$. Our main theorem then connects the bounds on the polyhedral norm on $\Re\widetilde\sigma_\Gamma$ to polyhedral bounds on the growth indicator function $\psi_\Gamma$ and the exponential decay rate of matrix coefficients of $L^2(\Gamma\backslash G)$.

\begin{theorem}
\label{thm:main}
Let $G$ be a real semisimple connected non-compact Lie group with finite center and $\Gamma < G$ a discrete subgroup. 
Then, for all $\mu\in \overline{\mathfrak{a}^\ast_+}$,
\begin{align}\label{main eq}
	\sup_{\lambda\in \widetilde \sigma_\Gamma} \|\Re \lambda\|_{\mathrm{poly},\mu} 
  = \theta_\Gamma(\mu) \geq \max(0,\delta_{\Gamma}'(\mu)).
\end{align}
and for all $\mu\in \mathfrak{a}^{\ast,\mathrm{Her}}\cap \overline{\mathfrak{a}^\ast_+}$:
\begin{align}\label{main eq2}
   \theta_\Gamma(\mu) =\max(0,
   \delta_{\Gamma}'(\mu)).
\end{align}
\end{theorem}
A particular choice for $\mu$ in the theorem is to take $\mu= \rho$. In this case, the modified critical exponent and the critical exponent are related by $\delta_{\Gamma}'(\rho) = \delta_{\Gamma}(\rho) -1$,
   Theorem~\ref{thm:main} then reads
  \[
   \sup_{\lambda\in \widetilde \sigma_\Gamma} \|\Re \lambda\|_{\mathrm{poly},\rho}
  = \max(\delta_{\Gamma, \rho} -1, 0).
  \]

\noindent We refer to Figure~\ref{fig:sl3} for a visualisation. Recall that a unitary representation is called \emph{tempered} if the matrix coefficients are in $L^{2+\vep}(G)$ for every $\vep>0$. It is well known that the temperedness of a unitary representation (whether or not it is tempered) is equivalent to certain decay properties of its matrix coefficients. Given two functionals, $\alpha, \beta$ on $\mathfrak{a}$, write $\alpha \le \beta$ if $\alpha(v)\le \beta(v)$ for all choices of $v\in \mathfrak{a}$. With this at hand, we conclude:
\begin{corollary}\label{cor:temperedness1}
 $L^2(\Gamma\backslash G)$ is tempered if and only if $\psi_\Gamma\leq\rho$.
\end{corollary}
This confirms a conjecture by Hee Oh and generalizes \cite[Theorem 1.6]{OhTemperedness} of Edwards and Oh.
They prove this result for the case of $\Gamma$ being Zariski-dense and the image of a Borel Anosov representation (i.e.~an Anosov representation w.r.t a minimal parabolic subgroup) and their proof is based on mixing results for Anosov subgroups by Edwards, Lee, and Oh \cite{ELO23}.
Corollary~\ref{cor:temperedness1} also extends work of Benoist and Kobayashi \cite{BK1} on tempered homogeneous spaces (see below).

The deviation from temperedness is usually measured by the property of a representation
being \emph{almost $L^p$} (see Section~\ref{sec:tempered}).
We refer to Proposition~\ref{prop:spectrum_and_decay} (iii) for a quantitative
statement connecting this deviation with $\theta_\Gamma(\rho)$
and therefore with the polyhedral bounds $\|\Re \sigma_\Gamma\|_{\mathrm{poly},\rho}$
which occur here naturally.

Let us denote by $\sigma(\Delta)$ the spectrum of the Laplace-Beltrami operator on $L^2(\Gamma\backslash G/K)$. In contrast to the rank one case, bounding the bottom of the Laplace spectrum a priori does not suffice in higher rank to obtain a  characterization of temperedness and non-temperedness of $L^2(\Gamma\backslash G)$, because in higher rank there are known examples of non-tempered representations that lead to Laplace eigenvalues bigger then $\|\rho\|^2$ (see e.g. \cite{Speh81} for $G=\SL_3(\R)$ and $G=\SL_4(\R)$).
However, based on  Theorem~\ref{thm:main}
we can prove that temperedness of $L^2(\Gamma\backslash G)$ is nevertheless equivalent to the bottom of the Laplace spectrum being $\|\rho\|^2$ and we obtain a refined version of Corollary~\ref{cor:temperedness1}:

\begin{corollary}\label{cor:temperedness}
Let $G$ be a real semisimple connected non-compact Lie group with finite center and $\Gamma < G$ a discrete subgroup, then the following statements are equivalent:
	\begin{enumerate}
		\item $\widetilde \sigma_\Gamma \subseteq i\mathfrak a^\ast$.
		\item 
	For all $\varepsilon>0$, there is $d_\varepsilon>0$ such that for all $f_1,f_2\in L^2(\Gamma\backslash G)^K$:
			\[
				|\langle (\exp v)f_1,f_2\rangle| \leq d_\varepsilon e^{\varepsilon \|v\|} e^{- \rho(v)} \|f_1\|_2\|f_2\|_2.
			\]
		\item $\psi_\Gamma\leq  \rho$.
		\item $L^2(\Gamma\backslash G)$ is almost $L^{2}$.
		\item $\min \sigma(\Delta)=\|\rho\|^2$.
		\item $L^2(\Gamma\backslash G)$ is tempered.
	\end{enumerate}
\end{corollary}
	In fact, (v) implies (iii) by \cite[Cor.~1.2]{WZ23}.
	More generally, in our notation
	\[
		\|\rho \|^2 - \max(0,\delta_{\Gamma}'(\mu))^2\left( \min_{H \in \overline{\mathfrak{a}_+}} \frac{\mu(H)}{\|H\|} \right)^2
		\geq \min \sigma(\Delta) \geq \|\rho\|^2 - \max(0,\delta_{\Gamma}'(\mu))^2 \|\mu\|^2
	\]
	for all $\mu\in \overline{\mathfrak{a}^\ast_+}$ \cite[Cor.~1.3]{WZ23}.

Note that, if $\Gamma$ is a lattice subgroup then, none of the above statements apply in this case and we get nothing novel from our result: $\psi_\Gamma = 2\rho$ and $L^2(\Gamma\bk G)$ always contains the trivial representation and is thus not tempered. Furthermore, the constant function leads to a zero eigenvalue of the Laplacian respectively a joint eigenvalue with spectral parameter $\rho$. For lattices one would instead have to study the temperedness of $L^2_0(\Gamma\backslash G)=\{f\in L^2(\Gamma\backslash G)| \int_{\Gamma\backslash G} f = 0\}$, but a general characterization of temperedness for this representation seems completely out of reach, given the fact that even in the special case of congruence subgroups of $\SL_2(\mathbb Z)$, the question of temperedness of $L^2_0(\Gamma\backslash  G)$ amounts to solving the longstanding Selberg conjecture \cite{Sel65, LRS95}.

As in the classical result of Patterson and Elstrodt, Theorem~\ref{thm:main}
provides a sharp bound on the size of the real part of the spectrum. In sharp
contrast to the rank one case, Theorem~\ref{thm:main}, however, provides bounds
on $\Re\wt \sigma_\Gamma$ with respect to different polyhedral norms and
relates them to different $\delta_\Gamma'(\mu)$. This is precisely the case when
$\dim(\mathfrak{a}^{*,\mathrm{Her}})>1$ or in other words if the root system of
the reduced restricted roots of $G$ is not of type $A_1$ or $A_2$. This fact can be
exploited to obtain the result below.
For its formulation we recall that the Benoist limit cone is
defined as
\[
	\mathcal{L}_\Gamma \coloneqq \left\{\lim_{i\to \infty} t_i \mu_+(\gamma_i) \in \overline{\mathfrak{a}_+} \mid
	t_i \to 0, \gamma_i\in \Gamma\right\}.
\]
\begin{theorem}\label{thm:limit_cone_tempered_intro}
	Let $G$ be of real rank $\geq 2$ not locally isomorphic to $\mathfrak{sl}_3(\mathbb{K}), \mathbb{K}=\R,\C,\mathbb{H},$ or $\mathfrak{e}_{6(-26)}$.
	Then, for every discrete subgroup $\Gamma< G$
satisfying $\mathcal L_\Gamma\subset \mathfrak{a}_+\cup \{0\}$, $L^2(\Gamma\backslash G)$ is tempered.
\end{theorem}
Combining this result with Corollary~\ref{cor:temperedness1} we deduce, that
\begin{equation}\label{eq:limit_cone_growth_indicator}
\mathcal L_\Gamma\subset \mathfrak{a}_+\cup \{0\} \Rightarrow \psi_\Gamma\leq \rho.
\end{equation}
This is an interesting implication because the limit cone a priori only encodes the directions in $\overline{\mathfrak a_+}$ in which there are infinitely many $\Gamma$-points, without saying anything about the growth rates of $\Gamma$-points in this direction. The latter is encoded in the growth indicator function. The implication~\eqref{eq:limit_cone_growth_indicator} thus says that whenever there are not infinitely many $\Gamma$-points in the direction of the walls of the Weyl chamber, the number of the $\Gamma$-points in the interior of the Weyl chambers can only grow with a moderate exponential rate.

A large class of subgroups fulfilling $\mathcal L_\Gamma\subset \mathfrak{a}_+\cup \{0\}$ is given by Borel Anosov subgroups, thus Theorem~\ref{thm:limit_cone_tempered_intro} proves the conjecture of Kim-Minsky-Oh \cite{KMO24}, except for the case of the $A_2$ root system for which the conjecture is still open. The conjecture of Kim, Minsky, and Oh was supported by the fact the they could prove the estimate $\psi_\Gamma\leq \rho$ for the class of Hitchin subgroups of $\SL(n,\R)$ based on the estimates on critical exponents by Potrie and Sambarino \cite{PS17}.
Note, however, that the assumption $\mathcal L_\Gamma\subset \mathfrak{a}_+\cup \{0\}$ is significantly weaker than being Borel Anosov and includes e.g.~all cusped Hitchin representations \cite{CZZ22}.
Furthermore note that with substantially more work, further conclusions on the interplay between the location of $\mathcal L_\Gamma$ and the shape of $\Re\widetilde\sigma_\Gamma$ as well as more properties of $\psi_\Gamma$ can be deduced with the help of Theorem~\ref{thm:main}.
A comprehensive study of these connections will be presented in \cite{Wol25}.

\subsection{Related Results}
\label{ss:rank one}

As discussed above, studying the connections between spectral properties of $\Gamma \bk G /K$ and the counting of $\Gamma$-points has a long history. The first instance of this connection is the characterization of the bottom $\inf\sigma(\Delta)$
of the Laplace spectrum for hyperbolic surfaces:
	\[ 
		\inf \sigma(\Delta) =\begin{cases}
			1/4&\colon\delta_\G<1/2\\
			1/4-(\delta_\G-1/2)^ 2&\colon \delta_\G\geq1/2,
		\end{cases} 
	\]
where $\delta_\Gamma$ is the critical exponent of the discrete subgroup $\Gamma < \SL_2(\R)$
\begin{equation}\label{eq:def_critical_exp}
	\delta_\G \coloneqq \inf\left\{s\in \R \colon \sum_{\g\in\G} e^{-s d(\gamma x_0, x_0)}<\infty \right\}, \quad x_0\in \mathbb{H}.
\end{equation}

This theorem is due to Elstrodt \cite{MR360472,MR360473,MR360474} and Patterson \cite{pattersonlimitset} and  has been extended to real hyperbolic manifolds of arbitrary dimension by Sullivan \cite{Sul87} and then to general locally symmetric spaces of rank one by Corlette \cite{Cor90}.

In our \emph{higher rank setting}, the bottom of the Laplace spectrum was estimated using the same definition of $\delta_\Gamma$ which is defined through $d(\gamma x_0,x_0)=\|\mu_+(x_0^{-1}\gamma x_0)\|$ by Leuzinger \cite{Leu04} and Weber \cite{Web08}. Later, Anker and Zhang \cite{AZ22} (see also \cite{CP04}) proved the exact formula
\[
	\inf \sigma(\Delta) = \begin{cases}
			\|\rho\|^2&\colon\tilde\delta_\G<\|\rho\|\\
			\|\rho\|^2-(\tilde\delta_\G-\|\rho\|)^ 2&\colon \tilde \delta_\G\geq\|\rho\|,
		\end{cases} 
	\]
	where  
 	$\tilde \delta_\Gamma$ is the modified critical exponent
	which is defined through $\|\mu_+(\gamma)\|$ and $\langle \rho,\mu_+(\gamma)\rangle$
	and therefore also takes the direction and not only the size of $\mu_+(\gamma)$ into account. However, as mentioned above, such bounds do not lead to temperedness of $L^2(\Gamma \bk G)$ due to the existence of non-tempered representations with arbitrary high Laplace eigenvalues.

A criterion of temperedness in higher rank was only achieved recently in the aforementioned work of Edwards and Oh \cite{OhTemperedness} for Borel Anosov representations.
Let us note that the generalization from Borel Anosov representations to general discrete subgroups is of great practical importance.
Already within the world of Anosov representations many concrete and important examples are Anosov with respect to a non-Borel parabolic subgroups e.g. holonomy groups of convex projective structures \cite{Benoist8, Benoist9} or maximal representations \cite{BIW}, see also \cite[Section 6]{GW}.
Also in some recent applications (\cite[Corollary 1.9]{DKO} or \cite{FO25}), that appeared since the present paper was in print, it was important to have Corollary~\ref{cor:temperedness1} for general Anosov subgroups.
Beyond the world of Anosov representations there are many classes of actively studied discrete subgroups such as relatively Anosov subgroups (see e.g. \cite{zz22}) or more generally $\theta$-transverse or $\theta$-divergence groups (see e.g. \cite{CZZ25} for a recent survey).

Temperedness in the complementary setting of homogeneous spaces $G/H$
for a closed subgroup $H$
with finitely many connected components has been studied by Benoist and Kobayashi in a series of papers \cite{BK1,BK2,BK3,BK4}.
They prove that the regular representation of $G$ on $L^2(G/H)$ is tempered if and only if a growth condition on $H$ is satisfied.
They also prove a version similar to Corollary~\ref{cor:temperedness}
(and also Proposition~\ref{prop:spectrum_and_decay})
where they
characterize when $L^2(G/H)$ is almost $L^p$ for $p\in 2\N$.

The main theorem (Theorem~\ref{thm:main}) not only gives a criterion for temperedness, but also allows one to locate the exceptional spectrum (i.e.~$\wt\sigma_\Gamma \cap (\mathfrak{a}_\C^\ast \setminus i \mathfrak{a}^\ast)$) via the sharp bounds on the polyhedral norms (see also Subsections \ref{ssec:sl3} and \ref{ssec:nontempered_exmpl} for an illustration in concrete examples). Such results were to our best knowledge not known in higher rank, except for the case where $G$ is a product of rank one groups and $\Gamma < G$ a general discrete,
torsion-free subgroup \cite{WW23}. The methods in \cite{WW23} however were based on analyzing the
resolvent kernels on the individual rank one factors and seemed not suitable for a generalization to general higher rank groups.

Concerning the quantitative bounds on the matrix coefficients, Kazhdan's Property (T) yields general estimates whenever $G$ has no factors locally isomorphic to $\mathfrak{so}(n,1)$ or $\mathfrak{su}(n,1)$. This amounts to a uniform bound on the quantities in \eqref{main eq},
i.e. an estimate independent of $\Gamma$, if $\Gamma$ has infinite covolume.
More precisely, in \cite[Thm.~7.1]{OhDichotomy} (see also previous work by Quint \cite{Qui03}) it is shown that $\psi_\Gamma\leq 2\rho-\Theta$
for some explicitly given functional $\Theta$.
Similarly, in \cite[Thm.~1.2]{OhUniformPointwise} it is shown that
\begin{align}
	\label{eq:decay propT}
	\left|\langle (\exp v)f_1,f_2\rangle_{L^2(\Gamma\backslash G)}\right|
	\leq C e^{-\Theta(v)}e^{\varepsilon\|v\|} \|f_1\|_2\|f_2\|_2
\end{align}
for all $ v\in \overline{\mathfrak{a}_+},$ and $ f_1,f_2 \in L^2(\Gamma\backslash G)^K$
for the same $\Theta$
(see also \cite{Li1995,LiZhu96}).
In \cite[Sect.~4A]{HWW23} one can find an analogues statement for the joint spectrum.
However, the bounds obtained by Property~(T) are not enough to deduce temperedness.
More precisely, the decay given by \eqref{eq:decay propT} is slower
than the decay required for temperedness,
as the functional $\Theta$ is in general smaller than $\rho$.
For example, $\Theta=\frac 12 \rho$ for $G=\SL_3(\R)$.
See \cite[Appendix]{OhUniformPointwise} for the precise values of $\Theta$ for the different root systems.

Let us finally mention two other recent results that concern the spectral theory of higher rank locally symmetric spaces of infinite volume: In \cite{EFLO23} Edwards, Fraczyk, Lee and Oh prove that the bottom of the Laplace spectrum is never an atom, provided that $\Gamma$ is a Zariski dense subgroup of infinite covolume in a semisimple real algebraic group $G$ with Kazhdan's property (T). They achieve this result by combining previous results on positivity of Laplace eigenvalues \cite{OhTemperedness} and the finiteness of Bowen Margulis Sullivan measures \cite{FL23}. In \cite{WW23b} the latter two named authors study the principal joint spectrum (i.e. the part of $\widetilde \sigma_\Gamma$ contained in $i\mathfrak a^*$) and give a dynamical criterion for the absence of embedded eigenvalues. Combining \cite[Theorem 1.1, Proposition 5.1]{WW23b} and Theorem~\ref{thm:main} we obtain:
\begin{corollary}
	Let $G$ be of real rank $\geq 2$ not locally isomorphic to $\mathfrak{sl}_3(\mathbb{K}), \mathbb{K}=\R,\C,\mathbb{H},$ or $\mathfrak{e}_{6(-26)}$.
 If $\Gamma$ the image of a Borel Anosov representation, then there exists no joint
 eigenfunction of the algebra of invariant differential operators $\mathbb
 D(G/K)$ in $L^2(\Gamma\bk G/K)$.
\end{corollary}

In the cases excluded in Theorem~\ref{thm:limit_cone_tempered_intro} and the previous
corollary,
 $\dim \mathfrak{a}^{\ast,\mathrm{Her}}=1$, i.e.~$\mathfrak{a}^{\ast,\mathrm{Her}}=\R\rho$.
In the these cases however, we can actually use this to locate where $\sup_{\lambda\in \wt \sigma_\Gamma} \|\Re \lambda\|_{\mathrm{poly},\rho}$ is attained.
\begin{proposition}
	\label{prop:sl3}  
	Let $G$ be locally isomorphic to $\mathfrak{sl}_3(\mathbb{K}), \mathbb{K}=\R,\C,\mathbb{H},$ or $\mathfrak{e}_{6(-26)}$
	and let $\Gamma < G$ be a discrete subgroup. Then the
	supremum $\sup_{\lambda\in \widetilde \sigma_\Gamma} \|\Re
	\lambda\|_{\mathrm{poly},\rho} = \max(0,\delta_{\Gamma}'(\rho))$ is achieved at
	$\lambda = \max(0,\delta_{\Gamma}'(\rho)) \rho$ (see Figure~\ref{fig:sl3}).
\end{proposition}
This in particular shows that there is a real spectral value on the boundary of the polyhedral region which is a priori not at all clear as we have no information on the imaginary part.

\subsection{Outline of the paper}
We start in Section~\ref{sec:preliminaries} with fixing the notation, introducing the joint spectrum of the algebra of invariant differential operators and recalling some important facts about tempered and almost $L^p$ representations. In Section~\ref{sec:decay_spec} we then study how the decay of matrix coefficients is related to the joint spectrum. A central tool is the Plancherel decomposition of $L^2(\Gamma\backslash G)$, the spherical dual, and asymptotics of spherical functions. The central result in this section is Proposition~\ref{prop:spectrum_and_decay}~(ii)
\[
\sup_{\lambda\in \widetilde \sigma_\Gamma} \|\Re \lambda\|_{\mathrm{poly},\mu}
= \theta_\Gamma(\mu), \quad \mu\in \overline{\mathfrak{a}_+^\ast},
\]
which roughly states that the polyhedral norm of the real part of the joint spectrum is determined by the decay of $L^2$-matrix coefficients. This relation should be known to experts, although we couldn't find any reference in the literature.

The central step of the paper is done in Section~\ref{sec:decay_and__growth_indicator}
where we derive a precise relation
between the decay of matrix coefficients for functions $f_1, f_2\in
C_c(\Gamma\backslash G)$ and the growth indicator function $\psi_\Gamma$
(Theorem~\ref{thm:coefficientdecay}). The decay of $C_c^\infty$-matrix coefficients without a uniform bound in $f_1,f_2$ is a priori significantly weaker than the uniform $L^2$-bounds described by $\theta_\Gamma(\mu)$. However, recent results of Cowling \cite{Cow23} allow us to pass from non-uniform $C_c$-bounds to the necessary uniform $L^2$-bounds. 
We conclude Section~\ref{sec:decay_and__growth_indicator} with the proof of Theorem~\ref{thm:main}.

In Section~\ref{sec:temperednesslimitcone} we prove Theorem~\ref{thm:limit_cone_tempered_intro} based on the existence of an optimal Hermitian functionl $\mu_\Gamma\in\mathfrak a^{*\mathrm{Her}}$ for the polyhedral estimates on Theorem~\ref{thm:main}.

Finally, in Section~\ref{sec:exmample}, we illustrate the implication of our main theorem for three concrete examples, the case of $G=\SL_3(\R)$, the product case, and recent example by Fraczyk and Oh \cite{FO25} in $\mathrm{SO}_0(2,n)$.

\emph{Acknowledgements:}
We thank Valentin Blomer for his suggestion to study this question and for numerous
stimulating discussions. We furthermore thank Michael Cowling, Samuel Edwards, Alex Gorodnik, Joachim Hilgert, Alex Kontorovich, Hee Oh, and Andres Sambarino for discussions and advice to the literature. Furthermore, we are very grateful to all of the referees for diligently reading our submission and providing informative feedback. This work has
received funding from the Deutsche Forschungsgemeinschaft (DFG) Grant No. SFB-TRR 358/1 2023 - 491392403 (CRC “Integral Structures in Geometry and Representation Theory”).

\section{Preliminaries}\label{sec:preliminaries}

\subsection{Notation}
In this article, $G$ is a real semisimple connected non-compact Lie group with finite center
and $K$ is a maximal compact subgroup of $G$, then $G/K$ is a Riemannian symmetric space of non-compact type. We fix an Iwasawa decomposition $G=KAN$, and have $A\cong \R^r$ where $r$ is the real rank of $G$ or the rank of the symmetric space $G/K$, respectively. Furthermore,  we define $M$ as the centralizer of $A$ in $K$ and $\overline{N}$ to be the nilpotent subgroup such that $KA\overline{N}$ is the opposite Iwasawa decomposition.
We denote by $\mathfrak{g},\mathfrak{k},\mathfrak{a},\mathfrak{n},\mathfrak{m},\overline{\mathfrak{n}}$  the corresponding Lie algebras.
For $g\in G$ let $H(g)\in \mathfrak{a}$ be the logarithm of the $A$-component in the Iwasawa decomposition.
Let $\Sigma \subseteq \mathfrak{a}^\ast$ be the root system of restricted roots,
$\Sigma^+$ the positive system corresponding to the Iwasawa decomposition,
and $W$ the corresponding Weyl group acting on $\mathfrak{a}^\ast$.
As usual, for $\alpha\in \Sigma$, we denote
by $m_\alpha$ the dimension of the root space,
and by $\rho$ the half sum of positive restricted roots counted with multiplicity.
Let $\mathfrak{a}_+ =\{H\in \mathfrak{a}\mid \alpha(H)>0 \:\forall \alpha\in \Sigma\}$ the positive Weyl chamber,
$\overline{\mathfrak{a}_+}$ its closure,
and $\mathfrak{a}_+^\ast$ the corresponding cone in $\mathfrak{a}^\ast$
via the identification $\mathfrak{a} \leftrightarrow \mathfrak{a}^\ast$ through the Killing form $\langle\cdot,\cdot\rangle$.
We have the Cartan decomposition $G=K\exp(\overline{\mathfrak a_+})K$
and for $g\in G$ there is a unique $\mu_+(g)\in \overline{\mathfrak a_+}$ such that
$g\in K \exp(\mu_+(g))K$.
For the Cartan decomposition the following integral formula holds (see \cite[Thm.~I.5.8]{gaga}):
\begin{equation}\label{eq:KAK_integral}
\int_{G}^{} {f(g)} \: dg = \int_{K}^{} {\int_{\mathfrak{a}_+}^{} {\int_{K}^{} {f(k\exp(H)k')\delta(H)} \: dk
} \: dH
} \: dk'
\end{equation}
where $\delta(H)=\prod_{\alpha\in \Sigma^+} (\sinh(\alpha(H))^{m_\alpha}$.
Note that $\delta(H)\leq e^{2\rho(H)}$.
We fix a discrete subgroup $\Gamma\leq G$.

\subsection{The growth indicator function}
\label{ss:gif}

In this subsection we recall the definition of the growth indicator function $\psi_\Gamma$.
It was introduced by Quint \cite{Qui02}
generalizing the critical exponent in rank at least $2$.
For an open cone $\mathcal{C}\subseteq \mathfrak{a}$,
let $\tau_{\mathcal{C}}$ be the abscissa of convergence for the series
$\sum_{\gamma\in \Gamma,\mu_+(\gamma)\in \mathcal{C}} e^{-s\|\mu_+(\gamma)\|}<\infty$,
	i.e. 
	\[
		\tau_{\mathcal{C}} = \inf\{s\in \R\mid 
		\sum_{\gamma\in \Gamma,\mu_+(\gamma)\in \mathcal{C}} e^{-s\|\mu_+(\gamma)\|}<\infty\}.
	\]
	The growth indicator function $\psi_\Gamma\colon \mathfrak{a}\to \R\cup \{-\infty\}$ is then defined as
	$\psi_\Gamma (u)= \|u\| \inf_{\mathcal{C}\ni u} \tau_{\mathcal{C}}$,
	where the infimum runs over all open cones $\mathcal{C}\subseteq \mathfrak{a}$
	containing $u$.
	We also set $\psi_\Gamma(0)=0$.
	One observes that $\psi_\Gamma$ is a positively homogeneous function
	that is upper semicontinuous.
	Moreover, $\psi_\Gamma$ is independent of the norm used on $\mathfrak{a}$.
	However, one usually uses the norm induced by the Killing form
	as it has the advantage of being invariant under the Weyl group
	which implies that $\psi_\Gamma$ is invariant under the opposition involution
	of $\mathfrak{a}$ given by $-w_0$,
	where $w_0$ is the Weyl group element with $w_0(\mathfrak{a}_+)=-\mathfrak{a}_+$.
	One also finds that $\psi_\Gamma\leq 2\rho$,
	$\psi_\Gamma=-\infty$ outside $\overline{\mathfrak{a}_+}$,
	and $\psi_\Gamma>-\infty$ implies $\psi_\Gamma \geq 0$.
The cone $\{v\in \mathfrak{a}\mid \psi_\Gamma(v)>-\infty\}$
	is precisely the limit cone
	\[
		\mathcal{L}_\Gamma = \{\lim_{i\to \infty} t_i \mu_+(\gamma_i) \in \overline{\mathfrak{a}_+} \mid
		t_i \to 0, \gamma_i\in \Gamma\}.
	\]
	If $\Gamma$ is Zariski-dense in a real algebraic group $G$,
	then one can make this more precise.
	Namely by \cite{Qui02}, $\psi_\Gamma>0$ on the interior of $\mathcal{L}_\Gamma$
	which is non-empty \cite{Ben96}
	and $\psi_\Gamma$ is concave.

\subsection{Algebra of invariant differential operators}
As mentioned in the introduction, $\mathbb{D}(G/K)$ denotes the algebra of $G$-invariant differential operators on $G/K$. The key result that allows a precise understanding of this algebra is the \emph{Harish-Chandra isomorphism} (see \cite[Thm.~II.5.18]{gaga}), for $\lambda \in \mathfrak{a}_{\C}^\ast$, let
\[
	\chi: \left\{\begin{array}{rcl}
			\mathbb{D}(G/K) & \longrightarrow &\mathrm{Poly}(\mathfrak{a}_\C^\ast)^W\\
D& \longmapsto & \{\lambda\mapsto \chi_\lambda(D), \lambda\in \mathfrak a_\C^\ast\}
	             \end{array}
	\right.
\]
which is an algebra isomorphism between $\mathbb D(G/K)$ and the algebra of Weyl group invariant polynomials on $\mathfrak a^*_\C$. In particular, one deduces that $\mathbb D(G/K)$ is abelian and is generated by $\mathrm{rank}(G/K)$ algebraically independent generators.

For any $\lambda\in\mathfrak a_\C^*$ we can define the \emph{elementary spherical function}
\begin{equation}\label{eq:spherical_function}
\phi_\lambda (g) \coloneqq \int_{K}^{} {e^{-(\lambda+\rho)H(g^{-1}k)}} \: dk,
\end{equation}
where 
$H\colon G\to \mathfrak a$ is defined by $g\in Ke^{H(g)}N$. This is a bi-$K$-invariant function and it descends to a left $K$-invariant function on $G/K$ which is a joint eigenfunction of $\mathbb D(G/K)$ fulfilling
\[
	D \phi_\lambda = \chi_\lambda(D) \phi_\lambda \quad \forall\: D\in\mathbb{D}(G/K).
\]
In fact, $\phi_\lambda$ is the unique such eigenfunction with $\phi_\lambda(e) = 1$ and for $\lambda,\lambda'\in\mathfrak a_\C^*$, $\phi_\lambda = \phi_{\lambda'}$ if and only if $\lambda '\in W\lambda$.

The elementary spherical functions $\phi_\lambda$ are parametrized by complex valued linear forms
$\lambda\in\mathfrak a^*_\mathbb C$ and, as the integral formula \eqref{eq:spherical_function} suggests, the imaginary part describes the oscillations of the spherical function, whereas the real part describes the asymptotic decay of its absolute value.
The latter can be explicitly expressed as follows \cite[Prop 7.15]{Knapp86}: For $\Re\lambda\in \overline{\mathfrak a^*_+}$ there are $C,d>0$ such that for all $H\in  \overline{\mathfrak a_+}$,
\[
 \left|\phi_\lambda(\exp(H))\right|\leq C e^{(\Re\lambda-\rho)(H)}(1+\rho(H))^d.
\]

Let us next study the action of $\mathbb D(G/K)$ on the locally symmetric space $\Gamma\bk G/K$: Each $D\in \mathbb{D}(G/K)$ is $G$-invariant
and therefore descends to $\Gamma\backslash G/K$.
All $D$ are unbounded operators on $L^2(\Gamma \backslash G/K)$, densely defined on $C_c^\infty(\Gamma\backslash G/K)$, and extend to normal operators on $L^2(\Gamma\backslash G/K)$ (we refer to \cite[Section 3.2]{WW23} for more details). Thus, we can define, for any $D$, its $L^2(\Gamma\backslash G/K)$-spectrum and denote it by $\sigma_{L^2}(D)\subset \C$. The spectral theory of $\mathbb{D}(G/K)$ is however, best described by a joint spectrum instead by the individual spectra and it is most convenient to parameterize this spectrum via the Harish-Chandra isomorphism by elements in $\mathfrak a^*_\C$:
\begin{definition}\label{def:joint_spectrum}
The joint spectrum of $\mathbb D(G/K)$ is defined by
\[
 \widetilde \sigma_\Gamma \coloneqq \{\lambda \in \mathfrak{a}_\C^\ast\mid \chi_\lambda(D) \in
	\sigma_{L^2}(D) \quad \forall D\in \mathbb{D}(G/K)\}  \subset \mathfrak a^*_\C.
\]
\end{definition}
In fact one can also choose a set of generators $D_1,\ldots, D_r$ of $\mathbb D(G/K)$, show that these are strongly commuting normal operators and consider their joint spectrum in the sense of \cite[Chapter 5]{Sch12}. This definition, however, coincides with the technically easier Definition~\ref{def:joint_spectrum} as shown in \cite[Proposition 3.6]{WW23}.

\subsection{Spherical dual and joint spectrum}\label{sec:spherical_dual}
Let us denote with $\widehat G$ the unitary dual of $G$, with $\widehat G_{\mathrm{sph}}\subset \widehat G$ the spherical dual of $G$, i.e.~the set of equivalence classes of irreducible unitary representations containing a non-zero $K$-invariant vector, and with $\widehat G_{\mathrm{tmp}}$ the tempered representations, i.e.~the support of the Plancherel measure of $L^2(G)$.

In the following we describe how $\widehat G_{\mathrm{sph}}$ can be parameterized by subset
of $\mathfrak{a}_\C^\ast/W$ (see \cite[Thm.~IV.3.7]{gaga}):
For $\pi\in \widehat G_{\mathrm{sph}}$ let $v_K$ be a normalized $K$-invariant vector.
Then the function $\phi\colon G\to \C, \phi(g)=\langle \pi(g)v_K,v_K\rangle$
is bi-$K$-invariant and positive definite,
i.e.~the matrix $(\phi(x_i^{-1}x_j))_{ij}$ is positive semidefinite for any choice of
finitely many $x_i\in G$.
Furthermore, $\phi$ is an eigenvector for each element
in the algebra $\mathbb D(G/K)$ of $G$-invariant differential operators
on $G/K$.

Therefore, $\phi=\phi_\lambda$ is an elementary spherical function
for $\lambda\in \mathfrak{a}_\C^\ast$.
Recall that $\phi_\lambda=\phi_\mu$ if and only if $W\lambda = W\mu$.
It can be shown that the mapping $\pi\mapsto W\lambda$ is a bijection
of $\widehat G_{\mathrm{sph}}$ onto the set
$\{\lambda \in \mathfrak{a}_\C^\ast /W \mid \phi_\lambda \text{ is positive definite}\}$.
We identify the two sets and write $\pi_\lambda$ for the representation corresponding to $\lambda\in \mathfrak{a}_\C^\ast/W$ with $\phi_\lambda$ positive definite.
In particular, for $\lambda\in \widehat G_{\mathrm{sph}}$  we have
$\langle \pi_\lambda(g) v,w\rangle = \phi_\lambda(g) \langle v,w\rangle$
if $v,w$ are both $K$-invariant.

Every positive definite function on $G$ is bounded by its value at $1$
and therefore $\widehat G_{\mathrm{sph}} \subseteq \operatorname{conv}(W\rho)+i\mathfrak{a}^\ast$ by \cite[Thm.~IV.8.1]{gaga}.
	Recall from the introduction that $\operatorname{conv}(W\mu)$ for $\mu\in \mathfrak{a}^\ast$
is the convex hull of the Weyl orbit $W\mu$ of $\mu$
which can be characterized by (see \cite[Lemma~IV.8.3]{gaga})
\[
	\operatorname{conv}(W\mu)= \{\|\lambda\|_{\mathrm{poly},\mu} \leq 1\}=\{\lambda\in \mathfrak{a}^\ast\mid \lambda(wH)\leq \mu(H) \:\forall H\in \mathfrak{a}_+, w\in W\}.
\]
Moreover, every positive definite elementary spherical function $\phi_\lambda$ 
is Hermitian,
i.e.~it satisfies
$\phi_\lambda(g^{-1}) = \overline{\phi_\lambda(g)}$.
As $\phi_\lambda(g^{-1})=\phi_{-\lambda}(g)$ and $\overline{\phi_\lambda(g) }=\phi_{\overline \lambda}(g)$,
we must have $W(-\lambda)=W\overline{\lambda}$.
Hence, $\widehat G_{\mathrm{sph}} \subseteq \{\lambda\in \mathfrak{a}_\C^\ast \mid \exists w\in W \colon w\lambda =-\overline{ \lambda}\} \eqqcolon \mathfrak{a}^{\ast,\mathrm{Her}}_\C$.
Furthermore, we define $\mathfrak{a}^{\ast,\mathrm{Her}}\coloneqq \mathfrak{a}^{\ast,\mathrm{Her}}_\C \cap \mathfrak{a}^\ast$.
The corresponding elementary spherical functions are Hermitian and have values in $\R_{>0}$.
They are used in \cite{Cow23} to bound matrix coefficients.

Let us now explain the relation of the joint spectrum of the invariant differential operators and the spherical dual: Consider the unitary representation $R$ on $L^2(\Gamma\backslash G)$ by right multiplication. By the abstract Plancherel theory, it can be decomposed into a direct integral of irreducible representations
\[
	(R,L^2(\Gamma \backslash G))\simeq \int_{X}^{\oplus} {\pi_x} \: d\mu(x),
\]
where $(X,\mu)$ is a measure space and
\[
 \pi:\left\{\begin{array}{rcl}
 X&\longrightarrow &\widehat G\\
 x&\longmapsto &\pi_x
 \end{array}
 \right.
\]
is a measurable map. We should think of $X$ as the Cartesian product of the unitary dual $\widehat G$ and a multiplicity space.
The joint spectrum of $\mathbb{D}(G/K)$ on $L^2(\Gamma \bk G/K)$ can now be expressed as follows:
\begin{proposition}
		[{\cite[Prop. 3.6]{WW23}}]
\label{prop:joint_spec}
	\begin{equation*}
		\widetilde \sigma_\Gamma = \supp (\pi_\ast \mu)\cap{\widehat G_{\mathrm{sph}}} \subseteq \widehat G_{\mathrm{sph}}\subset \mathfrak{a}_\C^\ast.
	\end{equation*}
\end{proposition}

\subsection{Temperedness and almost \texorpdfstring{$L^{p}$}{Lp}}\label{sec:tempered}
Recall that a unitary $G$-representation
$(\tau, \mathcal  H)$ with Plancherel decomposition
\[
 (\tau, \mathcal  H)\simeq \int_{X}^{\oplus} {\pi_x} \: d\mu(x)
\]
is called \emph{tempered} if $\supp (\pi_*\mu) \subset \widehat G_{\mathrm{tmp}}\subset \widehat G$. Temperedness of unitary representations has many equivalent characterizations and we recall those that are relevant for this paper:

\begin{definition}
Let $p\geq 2$. A unitary representation $(\tau,\mathcal  H)$ of $G$ is called \emph{strongly $L^{p+\varepsilon}$} or \emph{almost $L^p$} if there is a dense
subset $V\subset\mathcal  H$  such that for any $v,w\in V$,
the matrix coefficient $g\mapsto \langle \tau(g) v,w\rangle$ lies in $L^q(G)$ for all $q>p$.
\end{definition}
Note that if $\tau$ is strongly $L^{p+\varepsilon}$,
then $\tau$ is also strongly $L^{q+\varepsilon}$ for any $q\geq p$  since any matrix coefficients are bounded.

Let us furthermore introduce the Harish-Chandra function $\Xi(g)=\phi_0(g) = \int_{K}^{} {e^{-\rho(H(gk))}} \: dk$. It is well-known that $\Xi$ is a smooth bi-$K$-invariant function of $G$ with values in $(0,1]$.
Furthermore, there is a constant $C$ such that 
\begin{equation}\label{eq:Xi_bound}
	e^{-\rho(H)}\leq \Xi(e^H)\leq C (1+ |H|)^d e^{-\rho(H)}
\end{equation}
for $H\in \mathfrak a_+$. Here, $d$ is the number of positive reduced roots.
Note that by \eqref{eq:KAK_integral} this implies that $\Xi\in L^{2+\varepsilon}(G)$ for every $\varepsilon>0$ \cite[\S4.6]{GangolliVaradarajan}
\begin{proposition}
	[{\cite[Thm.~1 and 2]{CHH88}}]\label{prop:tempered_equivalences}
Let $(\tau, \mathcal  H)$ be a unitary $G$-representation 
then the following are equivalent
\begin{enumerate}
	\item $(\tau,\mathcal  H)$ is tempered.
	\item $(\tau,\mathcal  H)$ is almost $L^{2}$.
	\item For any $K$-finite unit vectors $v,w\in \mathcal  H$ ,
	\[
		|\langle \tau(g) v,w\rangle|\leq \left(\dim\langle Kv \rangle \dim \langle Kw\rangle\right)^{1/2}
		\Xi(g),
	\]
	for any $g\in G$, where $\langle Kv\rangle$ denotes the subspace spanned by $\tau(K)v$.
	\end{enumerate}
\end{proposition}
Note that in \cite{CHH88} the group $G$ is assumed to be a semisimple algebraic group
over a local field. 
However, as observed in \cite{Sun09} the same holds without any modification of the proof
as soon as $G$ admits an Iwasawa decomposition.
The same applies to Proposition~\ref{prop:pointwise_from_Lp_integer} below.

Since we are not only interested in temperedness, being strongly $L^{p+\varepsilon}$
gives us a measure for the extent of the non-tempered part.
However, the connection to uniform pointwise bounds seems to be established 
only for $p\in 2\N$:
\begin{proposition}
	[{\cite[Cor.~on p.~108]{CHH88}}]
	\label{prop:pointwise_from_Lp_integer}
	If $\tau$ is a unitary representation without a non-zero invariant vector
	that is strongly $L^{2k+\varepsilon}$, $k\in \N$, 
	then for any $K$-finite unit vectors $v$ and $w$,
	\[
		|\langle \tau(g) v,w\rangle|\leq \left(\dim\langle Kv \rangle \dim \langle Kw\rangle\right)^{1/2}
		\Xi^{1/k}(g).
	\]
\end{proposition}
Clearly, since $\Xi\in L^{2+\varepsilon}(G)$ the opposite implication holds as well.

\section{Decay of coefficients and the joint spectrum}\label{sec:decay_spec}
The aim of this section is to work out how the decay of matrix coefficients is linked to the joint spectrum. We will in particular show that	$L^2(\Gamma\backslash G)$ is tempered if and only if $\widetilde \sigma_\Gamma\subseteq i\mathfrak{a}^\ast$
and that there is a relation between polyhedral bounds on $\Re\widetilde\sigma_\Gamma$ and the decay of matrix coefficients of $L^2(\Gamma\backslash G)$. As tools we use standard representation theory and asymptotics of spherical functions. Although we assume these relations to be known to experts, we include the statements and proof in order to make the article self-contained.

We first prove that bounds on the real part of the joint spectrum lead to decay estimates for the matrix coefficients.
\begin{lemma}
	\label{la:matrix_bound_spectrum}
	For all $\varepsilon>0$, there is $d_\varepsilon>0$ such that for all $f,g\in L^2(\Gamma\backslash G)^K$ and all $v\in \overline{\mathfrak{a}_+}$ we have
	\[
		|\langle R (\exp v)f,g\rangle| \leq 
		d_\varepsilon e^{\sup_{\lambda\in\widetilde \sigma_\Gamma} (\Re \lambda-\rho)(v)} e^{\varepsilon \|v\|} \|f\|_2\|g\|_2.
	\]
\end{lemma}
\begin{proof}
	We decompose $f,g\in L^2(\Gamma\backslash G)^K$	
	into $\int_{X}^{\oplus} {f_x} \: d\mu(x)$ and $\int_{X}^{\oplus} {g_x} \: d\mu(x)$, respectively, 
	according to the decomposition $L^2(\Gamma\backslash G)\simeq \int_{X}^{\oplus} {\pi_x} \: d\mu(x)$.
	Since $f$ and $g$ are $K$-invariant, $f_x$ and $g_x$ are contained in
	$\pi_x^K$ for $\mu$-almost every $x\in X$ and hence they vanish for almost
	every $x\in X$ with $\pi_x\notin \widehat G_{\mathrm{sph}}$.
	We thus get
	\[
		\langle R (\exp v) f,g\rangle =\int_{X}^{} {\langle \pi_x(\exp v)f_x,g_x\rangle} \: d\mu(x)
		=\int_{\pi^{-1}(\widehat G_{\mathrm{sph}})}^{} {\langle\pi_x(\exp v)f_x,g_x\rangle} \: d\mu(x).
	\]
	We recall that if $\lambda\in \mathfrak{a}^\ast_\C/W$ corresponds to $\pi_\lambda\in \widehat G_{\mathrm{sph}}$ we have
	\[
		\langle \pi_\lambda (g) v_K, v_K\rangle = \phi_\lambda(g) \langle v_K,v_K\rangle,
	\]
	for all $v_K\in \pi_\lambda^K$.
	Therefore,
	\[
		\langle R(\exp v) f,g\rangle 
		=\int_{\pi^{-1}(\widehat G_{\mathrm{sph}})}^{} {\phi_{\lambda_x}(\exp v)\langle f_x,g_x\rangle} \: d\mu(x).
	\]
	Hence we can estimate
	\begin{align*}
		|\langle R(\exp v) f,g\rangle|
		&\leq \int_{\pi^{-1}(\widehat G_{\mathrm{sph}})}^{} {|\phi_{\lambda_x}(\exp v)|\|f_x\|\|g_x\|} \: d\mu(x)\\
		&\leq \operatorname{ess sup}_{\pi_\ast \mu|_{\widehat G_{\mathrm{sph}}}} |\phi_{\lambda_x}(\exp v)| \|f\|_2 \|g\|_2
		\\
		&\leq \sup _{\lambda\in \widetilde \sigma_\Gamma} |\phi_\lambda (\exp v)| \|f\|_2 \|g\|_2.
	\end{align*}
	For the elementary spherical function we have the well-known estimates \cite[Prop.~4.6.1]{GangolliVaradarajan} 
\[
	|	\phi_\lambda(\exp v)|\leq e^{\Re \lambda (v)} \Xi(\exp v)
	\leq d_\varepsilon e^{\Re \lambda (v) } e^{-\rho(v)} e^{\varepsilon\|v\|}
\]
for $\Re \lambda\in \overline{\mathfrak{a}_+^\ast}$ and any $\varepsilon >0$.
	This completes the proof.	
\end{proof}

We also prove an inverse statement that shows that decay of matrix coefficients in $L^2(\Gamma\bk G)$ implies the existence of obstructions on the joint spectrum.

\begin{lemma}
	\label{la:modified_bound_spectrum_from_matrix}
Suppose that there exists a homogeneous function $\theta\colon \mathfrak{a}_+\to \R$ such that for all $\varepsilon>0$, there is $d_\varepsilon>0$
such that for any $K$-invariant functions $f,g\in L^2(\Gamma\backslash G)$ and
any $v\in \mathfrak{a}_+$
\[
|\langle R(\exp v)f ,g\rangle |\leq d_\varepsilon e^{-\theta(v)}e^{\varepsilon \|v\|}\|f\|_2\|g\|_2.
\]
This then implies that
\[
			\Re \lambda \leq \rho -\theta,
\]
for all $\lambda\in\widetilde \sigma_\Gamma$.
\end{lemma}
\begin{proof}
	Let $\tilde \varepsilon >0$, $X_{\mathrm{sph}} = \pi^{-1} (\widehat G_{\mathrm{sph}})$,
	$\lambda_0\in \widetilde \sigma_\Gamma$,
	and $A_{\tilde \varepsilon}: = \{x\in X_{\mathrm{sph}} \mid |\lambda_x - \lambda_0| <\varepsilon\}$.
	Then $\mu(A_{\tilde \varepsilon})>0$ by Proposition~\ref{prop:joint_spec}.
	Put $f_{\tilde \varepsilon}=\mu(A_{\tilde \varepsilon})^{-1/2} \int_{X}^\oplus \mathbbm 1_{A_{\tilde \varepsilon}}(x) w^K_x\: d\mu(x)$
	where $w_x^K\in \pi_x^K$ is normalized.
	By definition $f_{\tilde \varepsilon}\in L^2(\Gamma \backslash G)^K$ is normalized
	and $\langle R(\exp v) f_{\tilde \varepsilon}, f_{\tilde \varepsilon}\rangle =
	\mu(A_\varepsilon)^{-1}
	\int_{A_\varepsilon}^{} {\phi_{\lambda_x}(\exp v)} \: d\mu(x)$.
	We infer that $\phi_{\lambda_0} (\exp v) =\lim _{{\tilde \varepsilon} \to 0}
	\langle R(\exp v) f_{\tilde \varepsilon}, f_{\tilde \varepsilon}\rangle$ and
	therefore, by the assumed bound on the matrix coefficients, we get
	$|\phi_{\lambda_0}(\exp v)| \leq d_\varepsilon e^{-\theta(v)} e^{\varepsilon \|v\|} $
	for any $\varepsilon >0$.
	Without loss of generality assume $\Re\lambda_0\in \overline{\mathfrak a_+^\ast}$.
	From \cite[Thm. 3.5 and proof of Thm. 10.1]{vdbS87}
	it follows that there is a polynomial $p(t)$ such that
	\[
		\phi_{\lambda_0}(\exp tv) p(t)^{-1} e^{-t(\lambda_{0}-\rho)(v)}
		\to 1 \quad \text{as} \quad t\to \infty.
	\]
	Hence, 
	\[ 
		1\leq \limsup_{t\to \infty} d_\varepsilon |p(t)|^{-1} e^{t(-\theta(v) + \varepsilon\|v\| - \Re \lambda_0(v) +\rho(v))},
		\]
		for any $\varepsilon >0$.
		We conclude
		\[
			-\theta(v) +\varepsilon\|v\|-\Re \lambda_0(v) +\rho(v) > 0
		\]
		and 
		\[
			\Re \lambda_0 \leq \rho -\theta.
		\]
		This completes the proof.
\end{proof}

In the next proposition we state how the polyhedral bounds on the spectrum are related to almost $L^p$ properties for $L^2(\Gamma\bk G)$.
We also obtain the equality of Theorem~\ref{thm:main} between the polyhedral norm of the spectrum 
and $\theta_\Gamma(\mu)$.
\begin{proposition}~
	\label{prop:spectrum_and_decay}
\begin{enumerate}
 \item $L^2(\Gamma\backslash G)$ is tempered if and only if $\widetilde \sigma_\Gamma \subseteq i\mathfrak{a}^\ast$.
 \item For all $\mu\in \overline{\mathfrak{a}_+^{\ast}}$, 
	 \[
		 \sup_{\lambda\in \widetilde \sigma_\Gamma} \|\Re \lambda\|_{\mathrm{poly},\mu} 
		 = \theta_\Gamma(\mu).
	 \]
  \item For $p_\Gamma\coloneqq \min\{p \geq 2\mid L^2(\Gamma\bk G) \text{ is almost } L^p\}$
	 we have:
\[
		 p_\Gamma \leq \frac{2}{1-\theta_\Gamma(\rho)}\leq 2 \lceil p_\Gamma/2 \rceil .
	 \]
\end{enumerate}
\end{proposition}
\begin{remark}
	We note that
	\[
		\theta_\Gamma(\mu) \inf _{H\in \mathfrak{a}_+} \frac{\mu(H)}{\rho(H)}\leq
		\theta_\Gamma(\rho) \leq \theta_\Gamma(\mu) \sup_{H\in \mathfrak{a}_+}	\frac{\mu(H)}{\rho(H)},
	\]
	for all $\mu\in \overline{\mathfrak{a}^{\ast}_+}$.
	Hence, one also obtains a statement on $p_\Gamma$
	by considering $\theta_\Gamma(\mu)$ instead of $\theta_\Gamma(\rho)$.
	However, it will not be sharp anymore even in the case $p_\Gamma\in 2\N$.
\end{remark}

\begin{proof}
	We start by proving (ii):
	By definition of $\theta_\Gamma(\mu)$ we have
	\[ 
		|	\langle R(\exp v) f_1,f_2\rangle |\leq d_{\epsilon,\varepsilon} e^{((\theta_\Gamma(\mu)+\epsilon) \mu-\rho)(v)+\vep \|v\|} \|f_1\|_2 \| f_2\|_2,
	\]
	for all $\epsilon, \varepsilon>0$, $v\in \overline{\mathfrak{a}_+}$, and $f_1,f_2\in L^2(\Gamma\backslash G)^K$.
	By Lemma~\ref{la:modified_bound_spectrum_from_matrix} this implies 
	\[
		\Re\lambda(v)\leq (\theta_\Gamma(\mu) +\epsilon) \mu(v),
	\]
	for every $v\in \overline{\mathfrak{a}_+}$ and $\lambda\in \widetilde \sigma_\Gamma$,
	i.e. $\|\Re\lambda\|_{\mathrm{poly},\mu}\leq \theta_\Gamma(\mu)$.
	On the other hand, we have 
	\[
		|\langle R (\exp v)f_1,f_2\rangle| \leq 
		d_\varepsilon e^{\sup_{\lambda\in\widetilde \sigma_\Gamma} (\Re \lambda-\rho)(v)} e^{\varepsilon \|v\|} \|f_1\|_2\|f_2\|_2
	\]
	by Lemma~\ref{la:matrix_bound_spectrum}.
	It follows that, if for some $\theta'\geq 0$ we have
	$\Re\lambda(v)\leq \theta'\mu(v)$
	for every $\lambda\in \widetilde \sigma_\Gamma$ and $v\in \overline{ \mathfrak{a}_+}$,
	then 
	$\theta_\Gamma(\mu)\leq \theta'$.
	We conclude
	\[
		\theta_\Gamma(\mu) = \inf \{\theta'\geq 0\mid
		\Re \lambda(v)\leq \theta' \mu(v)
	\:\forall\: v\in \overline{ \mathfrak{a}_+}, \lambda\in \widetilde \sigma_\Gamma\}
	= \sup_{\lambda\in \widetilde\sigma_\Gamma} \|\Re\lambda\|_{\mathrm{poly},\mu}.
	\]

	To prove (iii): Let $q>2/(1-\theta_\Gamma(\rho))$ and consider $f_1,f_2\in
C_c(\Gamma\bk G) \subset L^2(\Gamma\bk G)$ which is a dense subspace. Then by
setting $\tilde f_i(g) \coloneqq \max_{k\in K} |f_i(gk)|$ we get right $K$-invariant
functions and compute
\begin{align*}
 \int_G |\langle R(g) f_1,f_2\rangle|^q dg &\leq \int_G \langle R(g) \tilde f_1,\tilde f_2\rangle^q dg,\\
 &\leq \int_{\mathfrak a_+} \langle R(\exp(H)) \tilde f_1,\tilde f_2\rangle^q e^{2\rho(H)} dH.
\end{align*}
We use the definition of $\theta_\Gamma(\rho)$ to obtain
\begin{align*}
  \int_G |\langle R(g) f_1,f_2\rangle|^q dg   \leq d_\varepsilon \|\tilde f_1\|_2\|\tilde f_2\|_2 \int_{\mathfrak a_+} e^{(q(\theta_\Gamma(\rho) -1) +2)\rho(H) +q\varepsilon\|H\|)}dH.
\end{align*}
By our choice of $q$ this is integrable for $\varepsilon$ sufficiently small.
	Consequently, $L^2(\Gamma\bk G)$ is almost $L^q$
	so that $p_\Gamma \leq q$ for all $q> 2/(1-\theta_\Gamma(\rho))$.
	This proves the first inequality.

Conversely, if $L^2(\Gamma\bk G)$ is almost $L^{2k}$ then, by Proposition~\ref{prop:pointwise_from_Lp_integer}, we get that,  for any $f_1,f_2\in L^2(\Gamma\backslash G)^K$,
\[
 |\langle R(g) f_1,f_2\rangle| \leq  \|f_1\|_2\|f_2\|_2 (\Xi(g))^{\frac{1}{k}}
\]
and thus by \eqref{eq:Xi_bound}, for any $\varepsilon>0$,
\[
 |\langle R(\exp(v)) f_1,f_2\rangle| \leq  d_\varepsilon  e^{-\frac{1}{k}\rho(v)}e^{\varepsilon\|v\|}\|f_1\|_2\|f_2\|_2.
\]
Consequently, $\theta_\Gamma(\rho) \leq 1-1/k$ or equivalently $k\leq 1/(1-\theta_\Gamma(\rho))$.
Choosing $k$ as $\lceil p_\Gamma / 2 \rceil$ proves (iii).

Finally, (i) follows from (ii) and (iii)
because temperedness is equivalent to being almost $L^2$.
\end{proof}

\section{Decay of matrix coefficients and the growth indicator function}\label{sec:decay_and__growth_indicator}
In this section we study the connection between the decay of matrix coefficients and the growth indicator function.
We start with a slight modification of \cite[Prop. 7.3]{OhDichotomy}.

\begin{lemma}
	\label{la:modified_bound_psi_from_matrix}
	Suppose there exists a lower semicontinuous, homogeneous function $\theta\colon \overline{\mathfrak{a}_+}\to \R$ such that, for any $\varepsilon>0$ and $f,g\in C_c(\Gamma\backslash G)^K$, there is a $d_{\varepsilon, f,g}>0$
	such that, for any $v\in\overline{ \mathfrak{a}_+}$,
\begin{equation}
 \label{eq:decay_to_psi}
|\langle R(\exp v)f ,g\rangle |\leq d_{\varepsilon, f,g} e^{-\theta(v) +\varepsilon \|v\|}.
\end{equation}
Then this implies
\[
	\psi_\Gamma\leq 2\rho -\theta.
\]
\end{lemma}
\begin{proof}

If we have $\psi_\Gamma(u) =-\infty$ there is nothing to prove. 
Hence, we can assume $\psi_\Gamma(u)\geq 0$
and that $u \in \overline{\mathfrak{a}_+}$ is normalized.
 Fix an open (with respect to the relative topology of $\overline{\mathfrak{a}_+}$)
 cone $\mathcal C\subseteq \overline{ \mathfrak a_+}$
 containing $u$, and set $\mathcal C_T=\{v\in \mathcal C:
 \|v\|\leq T\}$ and $B_{T}=K\exp (\mathcal C_T) K$ for each $T>1$.
 
 Let $\varepsilon>0$.
 Let $U_\varepsilon=K \exp\{v\in \mathfrak{a}\mid \|v\|\leq \varepsilon/2\} K$.
 Then $U_\varepsilon$ is a symmetric
 open neighborhood of the identity in $G$
 which injects to $\Gamma\backslash G$ 
 for sufficiently small $\varepsilon$.

 Define $$F_{T,\varepsilon}(g, h)\coloneqq \sum_{\gamma\in \Gamma}\mathbbm 1_{U_\varepsilon B_T U_\varepsilon} (g^{-1}\gamma h)$$ which we regard as a function on $\Gamma\backslash G\times \Gamma\backslash G$.
  As $B_T \subseteq g U_\varepsilon B_T U_\varepsilon h^{-1}$ for all $g,h\in U_\varepsilon$,
  we have $\# \Gamma\cap B_T \leq F_{T,\varepsilon} (g,h)$ for all $g,h\in U_\varepsilon$.
 Let $\Phi_\varepsilon$ be a non-negative $K$-invariant continuous function supported in $\Gamma\backslash \Gamma U_\varepsilon$ with $\int_{\Gamma\backslash G} \Phi_\varepsilon d(\Gamma g)=1$.

We compute
\begin{align*}
	\#\Gamma\cap B_T    &   \leq \int_{\Gamma\backslash G \times \Gamma\backslash G}
     F_{T,\varepsilon} (\Gamma g, \Gamma h) \Phi_\varepsilon (\Gamma g)\Phi_\varepsilon (\Gamma h) d(\Gamma g) d(\Gamma h),
     \\
     &=\int_{\Gamma\backslash G\times\Gamma\backslash G}\sum_{\gamma\in\Gamma}\mathbbm{1}_{U_\varepsilon B_{T}U_\varepsilon}(g^{-1}\gamma h)\Phi_\varepsilon(\Gamma g)\Phi_\varepsilon(\Gamma h)\,d(\Gamma g)\,d(\Gamma h),
     \\
     &=\int_{\Gamma\backslash G}\int_{G}\mathbbm{1}_{U_\varepsilon B_{T} U_\varepsilon }(g^{-1} h),
     \Phi_\varepsilon(\Gamma g)\Phi_\varepsilon(\Gamma h)\,dg\,d(\Gamma h),
     \\
     &=\int_{\Gamma\backslash G}\int_{G}\mathbbm{1}_{U_\varepsilon B_{T} U_\varepsilon }(g^{-1} )
     \Phi_\varepsilon(\Gamma hg)\Phi_\varepsilon(\Gamma h)\,dg\,d(\Gamma h), \\
     &=\int_{U_\varepsilon B_T U_\varepsilon }\left(\int_{\Gamma\backslash G}\Phi_\varepsilon(\Gamma h g^{-1})\Phi_\varepsilon(\Gamma h)\, d(\Gamma h)\right)\, dg,\\
     &=\int_{U_\varepsilon B_T U_\varepsilon }\langle R(g^{-1}) \Phi_\varepsilon, \Phi_\varepsilon \rangle\, dg,
     \\
     &=\int_{U_\varepsilon B_T U_\varepsilon }\langle  \Phi_\varepsilon, R(g)\Phi_\varepsilon \rangle\, dg.
\end{align*}
The set $U_\varepsilon B_T U_\varepsilon$ is contained in $K \mathcal{C}_{T,\varepsilon} K$,
where $\mathcal{C}_{T,\varepsilon} = \{v\in \overline{\mathfrak{a}_+}\mid d(v,\mathcal{C}_T) <\varepsilon\}$.
Therefore, by \eqref{eq:KAK_integral},
\begin{align*}
     \#\Gamma\cap B_T &\leq 
\int_{\mathcal C_{T,\varepsilon}} \langle \exp v. \Phi_\varepsilon, \Phi_\varepsilon \rangle 
     e^{2\rho (v) }dv, \\
     &
     \leq d_{\varepsilon,\Phi_\varepsilon,\Phi_\varepsilon} \int_{\mathcal C_{T,\varepsilon}} 
     e^{(2\rho-\theta)  (v) +\varepsilon \|v\|}dv.
\end{align*}
Furthermore, given $\varepsilon>0$ and $\mathcal{C}$, there is a compact set $B\subseteq \overline{\mathfrak{a}_+}$ and an open cone $\mathcal{C}'\supseteq \mathcal{C}$ such that $\mathcal{C}_{T,\varepsilon} \subseteq B \cup \mathcal{C}_{T+\varepsilon}'$,
where $\mathcal{C}_{T}'$ is defined similarly to $\mathcal{C}_T$ replacing $\mathcal{C}$ by $\mathcal{C}'$.
We infer
\begin{align*}
     \#\Gamma\cap B_T &\leq 
     C_\varepsilon + d_{\varepsilon,\Phi_\varepsilon,\Phi_\varepsilon}d_\varepsilon \int_{\mathcal{C}_{T+\varepsilon}'}^{} {
     e^{(2\rho-\theta)  (v) +\varepsilon \|v\|}}\: dv, \\
		      &\leq C_\varepsilon +  d_{\varepsilon,\Phi_\varepsilon,\Phi_\varepsilon}
		      \int_0^{T+\varepsilon} \int_{v\in \mathcal C', \|v\|=1} 
		      e^{(2\rho- \theta)  (tv) +\varepsilon t \|v\|} t^{\dim \mathfrak{a}-1} dv dt, 
 \\
		      &\leq C_\varepsilon +d_{\varepsilon,\Phi_\varepsilon,\Phi_\varepsilon}
		      (T+\varepsilon)^{\dim \mathfrak{a}-1} 
		      \mathrm{vol}(\{v\in \mathcal{C}'\mid \|v\|=1\})
		      \int_{0}^{T+\varepsilon} {e^{t(\eta+\varepsilon)}} \: d{t}, 
\\
		      &\leq C_\varepsilon +d_{\varepsilon,\Phi_\varepsilon,\Phi_\varepsilon}(T+\varepsilon)^{\dim \mathfrak{a}-1} 
		      \mathrm{vol}(\{v\in \mathcal{C}'\mid \|v\|=1\})
		      \frac 1 {\eta+\varepsilon} (e^{(T+\varepsilon)(\eta+\varepsilon)} -1),     
 \end{align*}
 where $\eta = \sup \{2\rho(v)-\theta(v)\mid v\in \mathcal{C}', \|v\|=1\}$.
 Therefore
 $$\limsup_{T\to \infty} \frac{\log \# (\Gamma\cap B_T)}T \leq \limsup_T \frac {(T+\varepsilon)(\eta +\varepsilon)} T= \eta +\varepsilon.$$ 
 On the other hand,
 as $\psi_\Gamma(u)\geq 0$,
 $$\psi_\Gamma(u)= \inf_{u\in \mathcal C} \limsup_{T\to \infty} \frac{\log \# (\Gamma\cap K\exp (\mathcal C_T) K)}T,$$ 
 where the infimum is taken over all open cones $\mathcal C$ containing $u$.
 As $\mathcal{C}$ shrinks to the ray $\R_+u$,
 we can also take $\mathcal{C}'$ shrinking to $\R_+ u$,
 so that we get 
 \[
	 \psi_\Gamma(u) \leq \limsup_{v\to u} 2\rho(v)-\theta(v)
	 =2\rho(u)- \liminf_{v\to u} \theta(v).
 \]
 As we assumed that $\theta$ is lower semicontinuous,
 i.e.~$\liminf_{v\to u} \theta(v) \geq \theta(u)$,
 we obtain $\psi_\Gamma(u)\leq 2\rho(u)-\theta(u)$
 as desired.
\end{proof}

As a direct consequence of Lemmas~\ref{la:modified_bound_psi_from_matrix} and \ref{la:matrix_bound_spectrum} we get the following proposition.
\begin{proposition}
\label{prop:bound_psi_from_spectrum}
\[
\psi_\Gamma ( v)\leq \sup_{\lambda\in \widetilde \sigma_\Gamma} \Re \lambda(v) + \rho(v).
\]
\end{proposition}

\begin{proof}
As $\Re \widetilde \sigma_\Gamma\subset \mathfrak a^*$ is a bounded set we deduce that $\theta(v)=\inf_{\lambda\in \widetilde \sigma_\Gamma} -\Re \lambda(v) +\rho(v)$ is a continuous homogeneous function and we can apply Lemma~\ref{la:modified_bound_psi_from_matrix} using the spectral bounds on the matrix coefficients given by Lemma~\ref{la:matrix_bound_spectrum}.
\end{proof}

Note that this bound on the counting function is even a little bit more precise compared to the bounds stated in the main theorem, because the right hand side is not simply a dilation of $\mu$ but might be a more precise functional.

For the converse, we prove

\begin{proposition}
	\label{prop:psiestimate_implies_bddcoefficients}
	For all $\varepsilon > 0$ and $\mu\in \mathfrak{a}^{\ast,\mathrm{Her}}\cap \overline{\mathfrak{a}^\ast_+}$
	with $\delta_{\Gamma}'(\mu)<\infty$,
	there is $d_{\varepsilon}>0$ such that,
	for all $f_1,f_2 \in L^2(\Gamma\backslash G)^K$ and all $v\in \overline{\mathfrak{a}_+}$ we have
\[
	|\langle R(\exp v) f_1,f_2\rangle| \leq
	d_\varepsilon e^{\varepsilon \|v\|} e^{\max(0,\delta_{\Gamma}'(\mu)) \mu(v)-\rho(v)}
	\|f_1\|_2\|f_2\|_2.
\]
\end{proposition}
A key ingredient is the following decay of matrix coefficients for compactly supported functions.
\begin{theorem}
	\label{thm:coefficientdecay}
	Let $\theta\colon \overline{\mathfrak{a}_+}\to \R$ be positively homogeneous and
	continuous with $\psi_\Gamma\leq \theta$.
	Then, for all $f_1,f_2\in C_c(\Gamma\backslash G)$ and $\varepsilon >0$,
	there exists $C>0$ such that
	\[
		|\langle R(\exp H) f_1, f_2\rangle _{L^2(\Gamma\backslash G)}|
		\leq C e^{\varepsilon \|H\|} e^{\theta(H)-2\rho(H)}
	\]
	for all $H\in \overline{\mathfrak{a}_+}$.
\end{theorem}
\begin{remark}
	If $H\notin \mathcal{L}_\Gamma$ and therefore $\psi_\Gamma=-\infty$
	in a neighborhood of $H$,
	we even have
	$\langle R(\exp(tH))f_1,f_2\rangle_{L^2(\Gamma\bk G)} = 0$
	for $t$ large enough.
\end{remark}

\begin{remark}
 It should further be noted, that if $\psi_\Gamma \leq \rho$, then the exponent in
 Theorem~\ref{thm:coefficientdecay} is smaller than in
 Proposition~\ref{prop:psiestimate_implies_bddcoefficients}, where the decay is
 studied for $L^2$ functions. This is a well known phenomenon, for example, for
 decay estimates for geodesic flows on convex co-compact hyperbolic surfaces
 with $\delta_\Gamma <\frac{1}{2}$.
\end{remark}
To pass from Theorem~\ref{thm:coefficientdecay} to the uniform bounds in Proposition~\ref{prop:psiestimate_implies_bddcoefficients},
we use the following result of Cowling.

\begin{lemma}
	[{\cite[Lemma 3.5]{Cow23}}]
	\label{lem:cowling}
	Let $\mu\in \mathfrak{a}^{\ast, \mathrm{Her}}$ and $(\pi, \mathcal{H})$ a unitary representation of $G$.
	Then the following statements are equivalent:
	\begin{enumerate}
		\item There is a dense subspace $\mathcal{H}^0$ of $\mathcal{H}$,
			such that for all $\xi$ and $\eta$ in $\mathcal{H}^0$,
			there is a constant $C(\xi,\eta)$ such that
			\[
				\left(\int_{K}^{} {\int_{K}^{} {|\langle \pi(kxk') \xi,\eta\rangle|^2} \: dk
				} \: dk'\right)^{1/2} \leq
				C(\xi,\eta) \phi_{\mu}(x)\quad \forall x\in  G;
			\]
		\item For all $\xi$ and $\eta$ in $\mathcal{H}$,
\[
				\left(\int_{K}^{} {\int_{K}^{} {|\langle \pi(kxk') \xi,\eta\rangle|^2} \: dk
				} \: dk'\right)^{1/2} \leq
				\|\xi\|_{\mathcal{H}} \|\eta\|_{\mathcal{H}} \phi_{\mu}(x)\quad \forall x\in  G.
			\]
	\end{enumerate}
\end{lemma}

\begin{proof}
	[Proof of Proposition~\ref{prop:psiestimate_implies_bddcoefficients}
	from Theorem~\ref{thm:coefficientdecay}]
Based on
Theorem~\ref{thm:coefficientdecay} we will show that the matrix coefficients
for functions in $C_c(\Gamma\backslash G)$ satisfy (i) of Lemma~\ref{lem:cowling}.
Let
$f_1,f_2\in C_c(\Gamma\backslash G)$. Since \begin{equation}
	\label{eq:reduction_to_K_inv}
	\left|\int_{\Gamma\backslash G}^{} {f_1(\Gamma g h) f_2(\Gamma g)} \: d(\Gamma g)\right|
	\leq \int_{\Gamma \backslash G}^{} {\max_{k\in K}|f_1(\Gamma ghk)| \max_{k\in K} |f_2(\Gamma gk)| } \: d(\Gamma g),
\end{equation}
we can assume that $ f_i$ is non-negative and right $K$-invariant.

Let $\vartheta \coloneqq \max(0,\delta_{\Gamma}'(\mu))
=\inf\{t\geq 0 \mid t\mu > \psi_\Gamma -\rho\}$
be finite,
so that $\psi_\Gamma\leq  \vartheta \mu +\rho\eqqcolon \theta$.
The right hand side is continuous
so that
we can apply Theorem~\ref{thm:coefficientdecay}.
Hence, we have
\[
|\langle R(\exp H) f_1, f_2\rangle _{L^2(\Gamma\backslash G)}|
\leq C_{\varepsilon,f_1,f_2} e^{\varepsilon \|H\| +(\vartheta\mu-\rho)(H)}
\leq C_{\varepsilon,f_1,f_2} \phi_{\vartheta \mu +c\varepsilon\rho}(\exp(H)),
\]
with $c>0$ such that $\|H\|\leq c \rho(H)$
(see e.g. \cite[Thm.~2.5]{Cow23} -- note that the polynomial term therein is $\geq 1$).

Since $f_i$ are $K$-invariant and
the elementary spherical functions are bi-$K$-invariant,
we also have
\begin{align*}
  \left(\int_{K}^{} {\int_{K}^{} {|\langle R(kxk') f_1,f_2\rangle_{L^2(\Gamma\bk G)}|^2} \: dk
  } \: dk'\right)^{1/2} &=
  |\langle R(\exp \mu_+(x)) f_1, f_2\rangle _{L^2(\Gamma\backslash G)}|
  \\
  &\leq
  C_{\varepsilon,f_1,f_2} \phi_{\vartheta\mu +c\varepsilon\rho}(x)
\end{align*}
for all $x\in G$.
We now apply Lemma~\ref{lem:cowling} to obtain
\[
|	\langle R(\exp H) f_1,f_2\rangle_{L^2(\Gamma\bk G)}|\leq
\phi_{\vartheta\mu + c\varepsilon \rho}(\exp H) \|f_1\|_2 \|f_2\|_2,
\]
for all $f_1,f_2\in L^2(\Gamma\backslash G)^K$.
Since
\[
	\phi_{\vartheta\mu + c\varepsilon \rho}(\exp H) \leq
	C \left(\prod_{\alpha\in \Sigma^+} (1+ \alpha(H))\right) 
	e^{\vartheta \mu(H)+(\varepsilon -1) \rho(H)}
	\leq C e^{\varepsilon \| H\|} 
	e^{\vartheta \mu(H)+(c\varepsilon-1) \rho(H)}
\]
(see again \cite[Thm.~2.5]{Cow23}) 
	and $\rho(v) \leq \|\rho\| \|v\|$, the proposition follows.
\end{proof}

Theorem~\ref{thm:coefficientdecay} follows from a compactness argument and the following lemma.

\begin{lemma}
	\label{la:coefficientdecay}
	Let $f_1,f_2\in C_c(\Gamma\backslash G)$, $H_0\in \overline{\mathfrak a_+}$ normalized,
	and $s>\psi_\Gamma (H_0)$.
	Then there exists $\delta>0$ and $C>0$ such that
	\[
		|\langle R(\exp tH) f_1, f_2\rangle _{L^2(\Gamma\backslash G)}|
		\leq C e^{t(s-2\rho(H))},
	\]
	for all $t\geq 0$ and $H\in B_\delta(H_0)$ normalized.
\end{lemma}

\begin{proof}
	[Proof of Theorem~\ref{thm:coefficientdecay}
	from Lemma~\ref{la:coefficientdecay}]
Let us fix an arbitrary $\varepsilon>0$.

For any $H_0\in \overline{ \mathfrak{a}_+}$, we can find an $s_{H_0}$ such that
$\psi_\Gamma(H_0)<s_{H_0}< \theta(H_0) + \varepsilon$.

Then by Lemma~\ref{la:coefficientdecay} for any $H_0\in \overline{\mathfrak{a}_+}$ normalized, 
there is $\delta>0$ and $C>0$ such that
	\[
		|\langle R(\exp tH) f_1, f_2\rangle _{L^2(\Gamma\backslash G)}|
		\leq C e^{t(s_{H_0}-2\rho(H))},
	\]
	for all $t\geq 0$ and $H\in B_\delta(H_0)$.
	
	By shrinking $\delta$, we can assume that $s_{H_0}< \theta(H)+\varepsilon$
	for any $H\in \overline{B}_\delta(H_0)$.
	Therefore,
	\[
		|\langle R(\exp tH) f_1, f_2\rangle _{L^2(\Gamma\backslash G)}|
		\leq C  e^{t\left(s_{H_0}-2\rho(H)\right)}
		\leq C e^{\theta(tH) +t\varepsilon - 2\rho(tH)}
	\]
	for $t\geq 0$ and $H\in \overline B_\delta(H_0)$
	with a constant depending on $s_{H_0}$ and $\delta$.

	By compactness of the unit sphere in $\mathfrak a$,
	we only need finitely many $H_0^i$ in order to have
	\[
	\overline{\mathfrak{a}_+} \subseteq \bigcup_i \R_+ \cdot \tilde B_i \text{ where } \tilde B_i\coloneqq  B_{\delta_i}(H_0^i) \cap\{H\in\mathfrak a, \|H\| = 1\}.
	\]
	Thus, the constant can be chosen uniformly proving the theorem.
\end{proof}

Before proving Lemma~\ref{la:coefficientdecay} let us prove the following lemma
that is certainly known to experts 
	(compare e.g.~with \cite[Prop.~3.7]{BK1})
but might still be of independent interest.

Recall, that by Bruhat decomposition (see \cite[Prop.~I.5.21]{gaga})
that the mapping
\[
	(\overline n, m,a,n)\mapsto \overline n man \in G
\]
is a bijection of 
$\overline N \times M \times A\times N$
onto an open submanifold of $G$ whose complement has Haar measure $0$.
Moreover,
\[
 \int_{G}^{} {f(g)} \: dg = \int_{\overline N\times M\times A\times N}^{} {f(\overline n man)e^{2\rho(\log a)}} \: d\overline n \: dm \: da \: dn.
\]

\begin{lemma}
	\label{lem:conjugationdecay}
	Let $\varphi_1,\varphi_2\in C_c(G)$ with $\supp \varphi_i\subseteq \overline{N}MAN$.
	Then
	there is a constant $C=C_{\varphi_1,\varphi_2}$ such that for all $h\in A$
	\[
		\left|\int_{G}^{} {\varphi_1(h^{-1}gh)\varphi_2(g)} \: dg\right| \leq C e^{-2|\rho(\log h)|}.
	\]
\end{lemma}
\begin{proof}
	By the triangle inequality we can assume that $\varphi_i \geq 0$.
	Since $\supp \varphi_i \subseteq \overline{N}MAN$	
	there exist compact sets $C_{\overline N}\subseteq \overline N$, $C_A\subseteq A$,
	and $C_N\subseteq N$ with $\supp \varphi_i\subseteq C_{\overline N} M C_A C_N$.
	We thus have
	\begin{align*}
		c& \coloneqq c_{\varphi_1, \varphi_2, h} \coloneqq \int_{G}^{} {\varphi_1(h^{-1} g h) \varphi_2(g)} \: dg,\\
&= 
\int_{C_{\overline{N}}\times M \times C_A \times C_N}^{} {\varphi_1(h^{-1}\overline{n}manh)\varphi_2(\overline{n}man)} e^{2\rho(\log a)} \: d\overline{n} \: dm \: da \: dn,\\
&\leq 
\|\varphi_2\|_\infty 
\int_{C_{\overline{N}}\times M \times C_A \times C_N}^{} {\varphi_1(h^{-1}\overline{n}manh)} e^{2\rho(\log a)} \: d\overline{n} \: dm \: da \: dn.
\end{align*}
Since $M$ centralizes $A$ and $A$ is abelian
\[
	c\leq 
\|\varphi_2\|_\infty 
\int_{C_{\overline{N}}\times M \times C_A \times C_N}^{} {\varphi_1(h^{-1}\overline{n}hmah^{-1}nh)} e^{2\rho(\log a)} \: d\overline{n} \: dm \: da \: dn.
\]
Estimating $\varphi_1$ by its absolute value and using that $A$ normalizes both $N$ and $\overline{N}$ we get
\begin{align*}
	c&\leq \|\varphi_1\|_\infty \|\varphi_2\|_\infty \int_{M}^{} {} \: dm \int_{C_A}^{} {e^{2\rho(\log a)}} \: da
	\int_{C_{\overline{N}}\cap hC_{\overline{N}}h^{-1}}^{} {} \: d\overline{n}
	\int_{C_{N}\cap hC_{N}h^{-1}}^{} {} \: dn,\\
&\leq \|\varphi_1\|_\infty \|\varphi_2\|_\infty \int_{M}^{} {} \: dm \int_{C_A}^{} {e^{2\rho(\log a)}} \: da
	\int_{C_{N}}^{} {} \: dn
	\int_{hC_{\overline{N}}h^{-1}}^{} {} \: d\overline{n}.
\end{align*}
Since the Jacobian factor for the diffeomorphism  $\overline{n}\mapsto h \overline{n} h^{-1}$ of $\overline{N}$ is $\det \operatorname{Ad}(h)|{\overline{\mathfrak n}}= e^{-2\rho(\log h)}$
	we have
	\[
		\int_{hC_{\overline{N}} h^{-1}}^{} {} \: d\overline{n}
		= \int_{\overline{N}}^{} {1_{ C_{\overline{N}}}(h^{-1} \overline{n}h) } \: d\overline{n}
	=\int_{\overline{N}}^{} {1_{C_{\overline{N}}}} (\overline{n}) e^{-2\rho(\log h)} \: d\overline{n}
	=\int_{C_{\overline{N}}}^{} {} \: d\overline{n} \:e^{-2\rho(\log h)}.	
		\]
		We conclude
		\begin{align*}
			c_{\varphi_1, \varphi_2, h}
&\leq \|\varphi_1\|_\infty \|\varphi_2\|_\infty \int_{M}^{} {} \: dm \int_{C_A}^{} {e^{2\rho(\log a)}} \: da
	\int_{C_{N}}^{} {} \: dn
	\int_{C_{\overline{N}}}^{} {} \: d\overline{n} \:e^{-2\rho(\log h)},\\
	&= C_{\varphi_1,\varphi_2} e^{-2\rho(\log h)}.
\end{align*}
Switching the roles of $N$ and $\overline{N}$ in the argument
gives an estimate $c_{\varphi_1,\varphi_2,h} \leq C_{\varphi_1,\varphi_2} e^{2\rho(\log h)}$
proving the lemma.
\end{proof}

Let us now prove Lemma~\ref{la:coefficientdecay}.

\begin{proof}
	[Proof of Lemma~\ref{la:coefficientdecay}]
Let $f_1,f_2\in C_c(\Gamma\backslash G)$.
We can find $\tilde  f_i\in C_c(G)$
such that $ f_i(\Gamma g)=\sum_{\gamma\in \Gamma} \tilde  f_i(\gamma g)$.

We then have 
\begin{align}
	\notag
	\langle R(h)f_1,f_2\rangle _{L^2(\Gamma\backslash G)}&=
	\int_{\Gamma\backslash G}^{} {f_1(\Gamma g h) f_2(\Gamma g)} \: d\Gamma g
	=  \int_{G}^{} {\tilde f_1(gh) f_2(\Gamma g)} \: dg, \\
	&=\sum_{\gamma\in \Gamma} \int_{G}^{} {\tilde f_1(gh) \tilde f_2(\gamma g)} \: dg.
	\label{eq:unfolding}
\end{align}

For any $g\in G$ there is an open neighborhood $U_g$ of $g$ such that
	$U_g^{-1} U_g \subseteq \overline{N}MAN$ since $\overline{N}MAN$ is an open neighborhood of the identity element.
	Since $\supp \tilde  f_i$ is compact there are finitely many $g_k$ such that
	$\supp \tilde  f_i \subseteq \bigcup_k U_{g_k}$.
	There exists a partition of unity $\chi_k$ subordinate to $U_{g_k}$,
	i.e.
	$\chi_k\in C_c(G)$ with $\supp \chi_k \subseteq U_{g_k}$ and $\sum_k \chi _k(x)=1$ for all $x \in \supp \tilde  f_i$.
	We decompose $\tilde  f_i$ as $\sum_k \chi_k \tilde  f_i$ in \eqref{eq:unfolding}.
	This allows us to assume without loss of generality that $\supp \tilde  f_i$ is contained in some $U_g$, since we can estimate each of the finite summands individually.
	In particular, we can assume that $(\supp \tilde  f_i)^{-1}\supp \tilde f_i \subseteq \overline{N}MAN$.

Let $\gamma\in \Gamma$ such that 
$
\int_{G}^{} {\tilde f_1(gh) \tilde f_2(\gamma g)} \: dg \neq 0.
$
Then there is $g\in G$ with $gh\in \supp \tilde f_1 $ and $\gamma g\in \supp \tilde f_2 $.
Therefore, 
$\gamma\in (\supp \tilde f_2)g^{-1} \subseteq \supp  \tilde f_2 h (\supp \tilde f_1)^{-1}$.
Hence, there are $s_1$ and $s_2$ in $\supp \tilde f_1$ and $\supp\tilde  f_2$, respectively,
with $\gamma=s_2 h s_1^{-1}$.
By a change of variables
\begin{align*}
	\int_{G}^{} {\tilde f_1(gh) \tilde f_2(\gamma g)} \: dg &=
\int_{G}^{} {\tilde f_1(gh) \tilde f_2(s_2 h s_1^{-1} g)} \: dg
=\int_{G}^{} {\tilde f_1((hs_1^{-1})^{-1}gh) \tilde f_2(s_2 g)} \: dg, \\
								     &=\int_{G}^{} {\tilde f_1(s_1 h^{-1}gh) \tilde f_2(s_2 g)} \: dg .
\end{align*}
If we define $\varphi_i(g) \coloneqq \max_{s \in \supp \tilde  f_i}  |\tilde  f_i(sg)|$ we can estimate
\[
\left|	\int_{G}^{} {\tilde f_1(gh) \tilde f_2(\gamma g)} \: dg\right| \leq
\int_{G}^{} {\varphi_1( h^{-1}gh) \varphi_2( g)} \: dg.
\]
Hence we have 
\[
|\langle R(h) f_1,f_2\rangle|
\leq \# (\Gamma \cap (\supp \tilde f_2) h (\supp \tilde f_1)^{-1}
)\int_{G}^{} {\varphi_1(h^{-1} g h) \varphi_2(g)} \: dg.
\]
Note that if $\varphi_i(g)\neq 0$ then there is $s\in \supp\tilde  f_i$ such that $sg\in \supp \tilde  f_i$.
Hence, $\supp \varphi_i \subseteq (\supp \tilde  f_i)^{-1} \supp \tilde  f_i$ is compact
and contained in $ \overline{N}MAN$.
Therefore, by Lemma~\ref{lem:conjugationdecay}
\[\int_{G}^{} {\varphi_1(h^{-1}gh)\varphi_2(g)} \: dg \leq C e^{-2\rho(\log h)}.\]
Lemma~\ref{la:coefficientdecay} now follows from Lemma~\ref{la:compactperturbationcartanprojection} and Lemma~\ref{la:boundgammapointsgrowthindicator} below.
\end{proof}
		
\begin{lemma}
	[see {\cite[Prop. 5.1]{Ben96}}]
	\label{la:compactperturbationcartanprojection}
	For any compact set $C\subseteq G$ there exists a compact set $L\subseteq \mathfrak{a}$ 
	such that $\mu_+(CgC) \subseteq \mu_+(g) +L$.
\end{lemma}

\begin{lemma}
	\label{la:boundgammapointsgrowthindicator}
	For all $H_0\in \mathfrak{a}_+$ normalized,
	all $L\subseteq \mathfrak a$ compact, all $t$ large enough, and all $s>\psi_\Gamma(H_0)$ 
	there exists $\delta>0$ 
	and $C>0$ such that
	\[
		\# \{\gamma\in \Gamma \mid \mu_+(\gamma)\in tH + L\} \leq C e^{ts}
	\]
	for $H\in B_\delta(H_0)$ normalized.
\end{lemma}
\begin{proof}
	If $\psi_\Gamma(H_0)<s<0$ then $H_0$ is not in the limit cone and $\psi_\Gamma(H_0)=-\infty$.
	Moreover, there is an open cone containing $H_0$ that contains only finitely many $\Gamma$ points.
	In particular, $\{\gamma\in \Gamma\mid \mu_+(\gamma)\in tH +L\}$ is empty for 
	$H\in B_\delta(H_0)$ and $t=t_{H_0}$ large enough depending on $H_0$.

	We now assume $s\geq 0$.
By definition there exists an open cone $\mathcal{C}$ containing $H_0$
	such that
	\[
		\sum_{\gamma\in \Gamma, \mu_+(\gamma)\in \mathcal C} e^{-s\|\mu_+(\gamma)\|} <\infty.
	\]
	Therefore, there is $C>0$ such that
	\[
		\# \{ \gamma\mid \mu_+(\gamma)\in \mathcal{C}, \|\mu_+(\gamma)\| \leq t\}
		\leq C e^{ts}.
	        \]
	Note that, for every $\delta>0$ with $\overline {B_\delta(H_0)}\subseteq \mathcal C$, 
	there is $t_0>0$ such that
	$tH+L\subseteq \mathcal{C}$ for every $t\geq t_0$ and $H\in B_\delta(H_0)$.
	If we take $R>0$ is such that $L\subseteq B_R(0)$ then we can estimate for all $t\geq t_0$ and $H\in B_\delta(H_0)$ normalized
	\begin{align*}
		\# \{\gamma\mid \mu_+(\gamma)\in tH +L\} \leq
		\#\{\gamma\mid \mu_+(\gamma)\in \mathcal{C}, \|\mu_+(\gamma)\|\leq t+R\}
		\leq (Ce^{sR}) e^{ts}. 
		\qquad\qedhere
	\end{align*}
\end{proof}

 \subsection*{Proof of Theorem~\ref{thm:main}}
 The equality of
 $\sup_{\lambda\in \widetilde \sigma_\Gamma} \|\Re \lambda\|_{poly,\mu}$
 and $\theta_\Gamma(\mu)$ for all $\mu\in \overline{\mathfrak{a}_+^\ast}$
 in \eqref{main eq} 
 is Proposition~\ref{prop:spectrum_and_decay} (ii).
 Lemma~\ref{la:modified_bound_psi_from_matrix} shows $\psi_\Gamma \leq \rho + \theta_\Gamma(\mu) \mu$,
 so that $\delta_{\Gamma}'(\mu)\leq \theta_\Gamma(\mu)$ for all $\mu\in \overline{\mathfrak{a}_+^\ast}$.
 Finally, Proposition~\ref{prop:psiestimate_implies_bddcoefficients}
 shows $\theta_\Gamma(\mu)\leq \max(0,\delta_{\Gamma}'(\mu))$ for all $\mu\in \mathfrak{a}^{\ast,\mathrm{Her}}\cap \overline{\mathfrak{a}_+^\ast}$ with $\delta_{\Gamma}'(\mu)<\infty$.

 \section{Temperedness and the limit cone}
 \label{sec:temperednesslimitcone}
We now want to draw some important implications of our sharp polyhedral norm estimates.

	Recall that Theorem~\ref{thm:main} provides estimates for the optimal convex polyhedra $\mathrm{conv}(W\mu)$ for different $\mu\in\overline{\mathfrak a^*_+}$. One might thus ask how the intersection of all these polyhedra looks like.
\begin{proposition}
	\label{prop:defmugamma}
 Let
 \[
		C_\Gamma\coloneqq \bigcap_{\substack{\mu\in \overline{\mathfrak{a}_+^\ast} \\ \Re \wt \sigma_\Gamma \subseteq \operatorname{conv}(W\mu)}}
		\operatorname{conv}(W\mu).
\]
then there is a unique $\mu_\Gamma\in \mathfrak a^{*,\mathrm{Her}}\cap \overline{\mathfrak a_+^\ast}$ such that $C_\Gamma = \mathrm{conv}(W \mu_\Gamma)$
\end{proposition}
We first prove the following general lemma on the intersection of convex Weyl group invariant  polyhedra:
	\begin{lemma}
		\label{la:intesectionconv}
		Let $S\subseteq \overline{\mathfrak{a}_+^\ast}$ be any set.
		Then there exists a (unique) $\mu_S\in \overline{\mathfrak{a}_+^\ast}$
		such that
		\[
			\bigcap_{\mu\in S} \operatorname{conv}(W\mu) = \operatorname{conv}(W\mu_S).
		\]
	\end{lemma}
	\begin{proof}
		By \cite[Lemma~IV.8.3]{gaga}, $\operatorname{conv}(W\mu) \cap \overline{\mathfrak{a}_+^\ast}
		= \{\lambda\in \overline{\mathfrak{a}_+^\ast}\colon \lambda(H)\leq \mu(H)
		\:\forall H\in \mathfrak{a}_+\}$.
		Hence, $\mu_S$ is unique.
				Let $\alpha_1,\ldots,\alpha_r$ be the simple roots associated with $\Sigma^+$
		and $c_i\colon \mathfrak{a}^\ast\to \R$ the linear maps given by
		$\lambda=\sum_{i=1}^r c_i(\lambda) \alpha_i$.
		Then, $\operatorname{conv}(W\mu) \cap \overline{\mathfrak{a}_+^\ast}
		=\{\lambda\in \overline{\mathfrak{a}_+^\ast}\colon
		c_i(\lambda)\leq c_i(\mu) \: \forall i\}$.
		Therefore,
		\[
			\bigcap_{\mu\in S} \operatorname{conv}(W\mu) \cap \overline{\mathfrak{a}_+^\ast} =
			\{\lambda\in \overline{\mathfrak{a}_+^\ast}\colon
			c_i(\lambda)\leq c_i(\mu) \: \forall i, \mu \in S\}
			=
			\{\lambda\in \overline{\mathfrak{a}_+^\ast}\colon
			c_i(\lambda)\leq\inf_{\mu\in S} c_i(\mu) \: \forall i\}.
		\]
		This leads to the definition of $\mu_S \in \mathfrak{a}^\ast$ by demanding
		$c_i(\mu_S)=\inf _{\mu\in S} c_i(\mu)$,
		i.e.~$\mu_S = \sum_{i=1}^r \left( \inf_{\mu\in S} c_i(\mu)\right) \alpha_i$.
		The lemma is proved if we can show $\mu_S \in \overline{\mathfrak{a}_+^\ast}$,
		i.e.~$\langle\mu_S, \alpha_i\rangle \geq 0$.
		Indeed, for all $\mu\in S$ we have
		$\langle\mu,\alpha_i\rangle\geq 0$.
		Hence,
		\[
			c_i(\mu) \|\alpha_i\|^2 \geq
			\sum_{j\neq i} c_j(\mu) (-\langle \alpha_j,\alpha_i\rangle) \geq
			\sum_{j\neq i} \left(\inf _{\mu\in S} c_j(\mu)\right) (-\langle \alpha_j,\alpha_i\rangle)
		\]
		as $\langle \alpha_j,\alpha_i\rangle \leq 0$ for $j\neq i$.
		This implies $
		\inf_{\mu\in S}c_i(\mu) \|\alpha_i\|^2 \geq
		\sum_{j\neq i} \left(\inf _{\mu\in S} c_j(\mu)\right) (-\langle \alpha_j,\alpha_i\rangle)
		$
		proving $\mu_S\in \overline{\mathfrak{a}_+^\ast }$.
	\end{proof}

	\begin{proof}[Proof of Proposition~\ref{prop:defmugamma}]
Lemma~\ref{la:intesectionconv} implies that there is a unique $\mu_\Gamma\in \overline{\mathfrak{a}_+^\ast}$ such that $C_\Gamma=\operatorname{conv}(W\mu_\Gamma)$.
	As $\Re \wt \sigma_\Gamma$ is invariant under opposition involution $\iota$ this is true for $C_\Gamma$.
	Hence, $\iota\mu_\Gamma = \mu_\Gamma$ and
	$\mu_\Gamma\in \mathfrak{a}^{\ast, \mathrm{Her}}$.
\end{proof}
Before proving Theorem~\ref{thm:limit_cone_tempered_intro} let us discuss its assumption on $G$.
\begin{remark}
	In the proof of Theorem~\ref{thm:limit_cone_tempered_intro}
	we will need $\dim \mathfrak{a}^{\ast,\mathrm{Her}}\geq 2$.
	The classification of irreducible root systems shows that this is the case
	if and only if the rank of the root system of reduced roots is $\geq 2$ and
	not of type $A_2$
	(see \cite[Plates~I-IX]{BourbakiLie} or \cite[Remark~4.4]{HWW23} where $\dim \mathfrak{a}^{\ast,\mathrm{Her}}=d_-$).
	From \cite[Appendix C]{Knapp86} we can read off that this root system only
occurs for the real Lie algebras $\mathfrak{sl}_3(\mathbb{K})$, $\mathbb{K}=\R,\C,\mathbb{H}$, and for E IV.
The latter is also denoted by $\mathfrak{e}_{6(-26)}$.
\end{remark}

\begin{proof}[Proof of Theorem~\ref{thm:limit_cone_tempered_intro}]
    Let us assume that $\mu_\Gamma\neq 0$.
    Then, by definition of $C_\Gamma$, it follows that
	\[
		\sup_{\lambda\in \wt \sigma_\Gamma}
		\|\Re \lambda\|_{\mathrm{poly},\mu_\Gamma} =
		\inf\{\theta\geq 0 \colon \Re \wt \sigma_\Gamma \subseteq \theta \operatorname{conv}(W\mu_\Gamma)\} =1
	\]
	and by Theorem~\ref{thm:main} that $\delta_\Gamma'(\mu_\Gamma) = 1$.
	Moreover, for each $\mu\in \overline{\mathfrak{a}^\ast_+}$,
	\[
		\operatorname{conv}(W\mu_\Gamma) \subseteq \sup_{\lambda\in \wt
		\sigma_\Gamma} \|\Re \lambda\|_{\mathrm{poly},\mu}\operatorname{conv}(W \mu),
	\]
	if the supremum is finite.
	Hence, without restriction on $\mu\in \overline{\mathfrak{a}_+^\ast}$,
	$\mu_\Gamma(v)\leq \sup_{\lambda\in \wt
	\sigma_\Gamma} \|\Re \lambda\|_{\mathrm{poly},\mu} \mu(v)$ for all $v\in \overline{\mathfrak{a}_+}$.
	In particular,
	\begin{equation}
		\label{eq:inequality_polynorm}
		\sup_{v\in \overline{\mathfrak{a}_+}} \frac{\mu_\Gamma(v)}{\mu(v)}
		\leq \sup_{\lambda\in \wt \sigma_\Gamma} \|\Re \lambda\|_{\mathrm{poly},\mu}.
	\end{equation}
	For $\delta_{\Gamma}'(\mu)$ we have:
	\begin{equation}
		\label{eq:etalimitcone}
		\delta_{\Gamma}'(\mu) = \sup_{v\in \overline{\mathfrak{a}_+}} \frac{\psi_\Gamma(v)-\rho(v)}{\mu(v)}
		= \sup_{v\in \mathcal{L}_\Gamma} \frac{\psi_\Gamma(v)-\rho(v)}{\mu(v)}
		\leq \sup_{v\in \mathcal{L}_\Gamma} \frac{\mu_\Gamma(v)}{\mu(v)}
	\end{equation}
	since $\delta_\Gamma'(\mu_\Gamma)=1$.
	Theorem~\ref{thm:main} implies
	$\sup_{\lambda\in \wt \sigma_\Gamma} \|\Re \lambda\|_{\mathrm{poly},\mu}
	= 	\delta_{\Gamma}'(\mu)$
	for all $\mu \in \mathfrak{a}^{\ast,\mathrm{Her}}\cap \overline{\mathfrak{a}^\ast_+}$
	as the left hand side is positive by the assumption $\mu_\Gamma\neq 0$.
	Combining this with \eqref{eq:inequality_polynorm} and \eqref{eq:etalimitcone}
	we get
	\begin{equation}
		\label{eq:equalitysups}
		\sup_{v\in \overline{\mathfrak{a}_+}} \frac{\mu_\Gamma(v)}{\mu(v)}
		=\sup_{v\in \mathcal{L}_\Gamma} \frac{\mu_\Gamma(v)}{\mu(v)}
		\quad \text{for all} \quad \mu\in
		\mathfrak{a}^{\ast,\mathrm{Her}}\cap \overline{\mathfrak{a}^\ast_+}.
	\end{equation}

	Since $\dim \mathfrak{a}^{\ast,\mathrm{Her}}\geq 2$
	we pick $\mu \in \mathfrak{a}^{\ast,\mathrm{Her}}\cap \overline{\mathfrak{a}^\ast_+}$
		such that $\mu\notin \R\mu_\Gamma$.
		Since $\mathcal{L}_\Gamma \subseteq \mathfrak{a}_+\cup \{0\}$
		and $\mu(w)>0$ for $w\in \mathfrak{a}_+$,
		$\mathcal{L}_\Gamma \cap \{\mu=1\}$ is compact.
		It follows that $\sup_{v\in \mathcal{L}_\Gamma} \frac{\mu_\Gamma(v)}{\mu(v)}
			=\mu_\Gamma(v_0)$ for some $v_0\in \mathcal{L}_\Gamma$ with $\mu(v_0)=1$.
		By \eqref{eq:equalitysups},
		\[
			\sup_{v\in \overline{\mathfrak{a}_+}} \frac{\mu_\Gamma(v)}{\mu(v)}
			=\sup_{\substack{v\in \overline{\mathfrak{a}_+}\\ \mu(v)=1}} \mu_\Gamma(v) 
			=\mu_\Gamma(v_0).
		\]
		It follows that $\mu_\Gamma$ vanishes on $\ker \mu$,
		i.e.~$\mu_\Gamma\in \R \mu$.
		This is a contradiction as we assumed $\mu_\Gamma\neq 0$ and $\mu \notin \R \mu_\Gamma$.
		Hence, we completed the proof of Theorem~\ref{thm:limit_cone_tempered_intro}
		since $\mu_\Gamma=0$ is equivalent to $\widetilde \sigma_\Gamma \subseteq i \mathfrak{a}^{\ast}$ and thereby to the temperedness of $L^2(\Gamma\bk G)$.
	\end{proof}

\section{Examples of precise descriptions of the spectrum}\label{sec:exmample}
In this last section we want to consider three concrete examples:
The product case $G=G_1\times G_2$ of two rank one groups, the case $G=\SL_3(\R)$, as well as a recently constructed concrete example of a non-tempered subgroup in $\mathrm{SO}(2,n)$.
In the product case we also consider the product of two discrete subgroups $\Gamma= \Gamma_1\times\Gamma_2$, such that the spectral theory of the joint spectrum of invariant differential operators trivially reduces to the rank one case.
Nevertheless we think that it is quite instructive to illustrate the main result in this concrete example.
In the case of $\SL_3(\R)$ we show that using the additional information of the root system $A_2$ with our main result allows us to deduce some finer information about the spectrum.
Additionally, the concrete non-tempered example in $\mathrm{SO}(2,n)$ provides a nice illustration how the polyhedral bounds on the spectrum yield quite precise information about the spectrum.

\subsection{Product case}
\label{sec:product_case}
Let us first consider the product case,
in which the joint spectrum is explicitly given
by the product of the two rank-one spectra
and which yields a nice illustration of our result:
More precisely, let $G=G_1\times G_2$ be the product of two rank one groups $G_i$, $i=1,2$.
We indicate by the subscript $i$ the corresponding subgroups of $G_i$ and the respective subspaces of their Lie algebras.
Assume that the discrete subgroup $\Gamma$ is also a product of discrete subgroups $\Gamma_i$ of $G_i$.
Clearly,
\begin{equation}\label{eq:product_spectrum}
 \widetilde \sigma_\Gamma = \big\{(\lambda_1,\lambda_2)\in \mathfrak{a}_{1,\C} \times \mathfrak{a}_{2,\C} \mid
		|\rho_i|^2-|\Re \lambda_i|^2 + |\Im \lambda_i|^2
	\in \sigma(\Delta_i)\big\}
	=\wt \sigma_{\Gamma_1}\times \wt \sigma_{\Gamma_2},
\end{equation}
\noindent where $\Delta_i$ is the Laplacian of $\Gamma_i\backslash G_i/K_i$ acting on one factor of $\Gamma\backslash G/K$.
Recall that $\inf \sigma(\Delta_i)= |\rho_i|^2 - \max(0,\delta_{\Gamma_i}-|\rho_i|)^2$  ,
where $\delta_{\Gamma_i}$ is the critical exponent of $\Gamma_i$.

\begin{figure}[ht]
\centering
\includegraphics[width=\textwidth, trim= 3cm 4cm 2cm 1cm, clip]{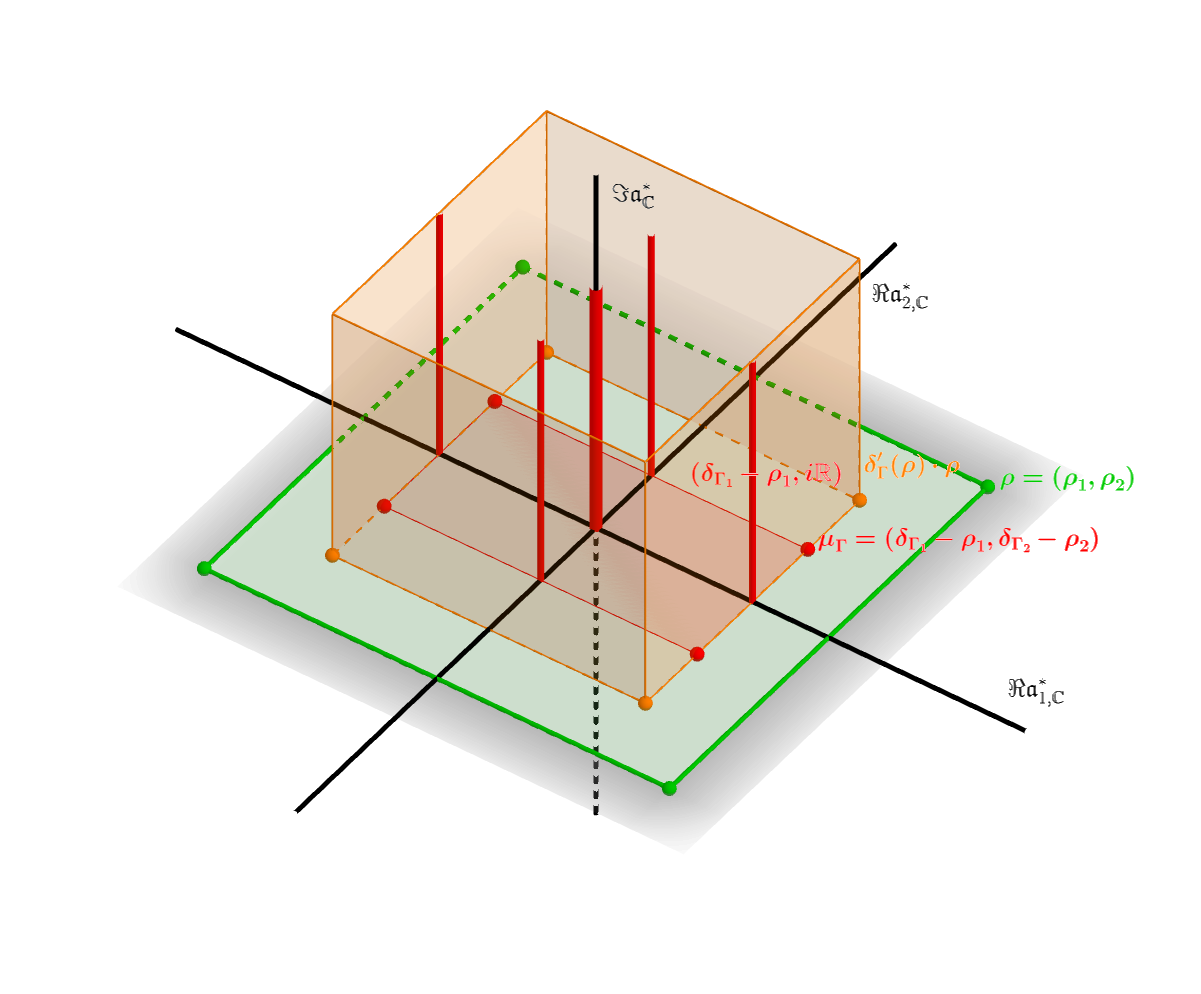}
\caption{Joint spectrum in the product case under the assumption that the single factors have no
	exceptional spectrum besides $\delta_{\Gamma_i}-\rho_i> 0$.
There is a joint eigenvalue $(\delta_{\Gamma_1}-\rho_1,\delta_{\Gamma_2}-\rho_2)$ but also continuous spectrum $\pm (\delta_{\Gamma_1}-\rho_1,i\R)$ and $\pm (i\R,\delta_{\Gamma_2}-\rho_2)$ 
as well as $i\R\times i\R$ (red).
One observes that the polyhedral bound with respect to $\rho=(\rho_1,\rho_2)$
gives the orange square, which is the smallest square containing $\Re \wt\sigma_\Gamma$.
However, the optimal rectangle is provided by considering $\mu_\Gamma=(\delta_{\Gamma_1}-\rho_1,\delta_{\Gamma_2}-\rho_2)$.
}
\label{fig:product}
\end{figure}

	\begin{lemma}
		\label{la:gifproduct}
	For $(H_1,H_2)\in\mathfrak{a}$, we have
	\[
		\psi_\Gamma(H_1,H_2)\leq \delta_1 |H_1| + \delta_2 |H_2|,
	\]
	and furthermore $\psi_\Gamma(H_1,0)=\delta_{\Gamma_1} |H_1|$ for $H_1\in \mathfrak{a}_{1,+}$,
	and $\psi_\Gamma(0,H_2)=\delta_{\Gamma_2} |H_2|$ for $H_2\in \mathfrak{a}_{2,+}$.
	Furthermore,
	if $\Gamma$ is Zariski-dense, then $
		\psi_\Gamma(H_1,H_2) = \delta_1 |H_1| + \delta_2 |H_2|,
		$
		for $(H_1,H_2)\in \overline{\mathfrak{a}_+}$.
\end{lemma}
\begin{proof}
	Let $(H_1,H_2)\in \overline{\mathfrak{a}_+}$ and assume $1 = |H_1|\geq |H_2|$.
	The opposite case is handled the same way.
	Let $c>0$ and $\mathcal{C} \subseteq \mathfrak{a}$ be the cone
	$\{(v_1,v_2)\in \overline{\mathfrak{a}_+} \colon ||v_2|/|v_1| -| H_2|| <c\}$.
	Then,	
	\begin{align*}
		\sum_{\mu_+(\gamma)\in \mathcal{C}} e^{-s\|\mu_+(\gamma)\|}
&\leq\sum_{\substack{(\gamma_1,\gamma_2)\in \Gamma \\
|\gamma_1|(|H_2|-c)<|\gamma_2| < |\gamma_1|(|H_2|+c)}} 
		e^{-s|\gamma_1| (1+(|H_2|-c)^2)^{1/2}},
		\\
	&	\leq \sum_{\gamma_1\in \Gamma_1} \#\{ \gamma_2\in \Gamma_2\mid 
		|\gamma_2| < |\gamma_1|(|H_2|+c)\}
		e^{-s|\gamma_1| (1+(|H_2|-c)^2)^{1/2}}.
	\end{align*}
	As $\#\{\gamma_2 \in \Gamma_2 \colon |\gamma_2|<R\}\leq 
	e^{(\delta_{\Gamma_2} +\varepsilon)R}$
	for $R$ big enough,
	this is finite if 
	\[
		\sum_{\gamma_1\in \Gamma_1} e^{(\delta_{\Gamma_2}+\varepsilon)|\gamma_1|(|H_2|+c)}
		e^{-s|\gamma_1| (1+(|H_2|-c)^2)^{1/2}}
		<\infty.
	\]
	This is the case if $s(1+(|H_2|-c)^2)^{1/2}-(\delta_{\Gamma_2} +\varepsilon)(|H_2|+c) > \delta_{\Gamma_1}$.
	Letting $\varepsilon\to 0$ and $c \to 0$ shows
$\psi_\Gamma(1,H_2) \leq(\delta_{\Gamma_1} +\delta_{\Gamma_2}|H_2|)$.

Conversely, for each cone $\mathcal{C}$
containing $(H_1,0)\in \overline{\mathfrak{a}_+}\setminus \{0\}$,
we have
\[
	\sum_{\gamma_1\in \Gamma_1} e^{-s|\gamma_1|} \leq \sum_{\mu_+(\gamma)\in \mathcal{C}}
e^{-s\|\mu_+(\gamma)\|}.
\]
Therefore, $\delta_{\Gamma_1} \leq \psi_\Gamma(H_1,0) |H_1|^{-1}$.
Hence, we have $\psi_\Gamma(H_1,0)=\delta_{\Gamma_1} |H_1|$.
Now, if $\Gamma$ is Zariski-dense and therefore $\psi_\Gamma$ is concave,
\[
	\psi_\Gamma(H_1,H_2)\geq 
	|H_1| \psi_\Gamma(H_1,0)+ |H_2| \psi_\Gamma(0,H_2)=
	\delta_{\Gamma_1}|H_1| + \delta_{\Gamma_2}|H_2|
\]
proving the lemma.
\end{proof}

By \eqref{eq:product_spectrum} it is clear that $\mu_\Gamma$ defined in Proposition~\ref{prop:defmugamma}
is given by $\mu_\Gamma(H_1,H_2) = \max(0,\delta_{\Gamma_1} - \rho_1)  |H_1| + \max(0,\delta_{\Gamma_2}-\rho_2) |H_2|$
for $(H_1,H_2)\in \overline{\mathfrak{a}_+}$.
If $\delta_{\Gamma_i} \leq \rho_i$ for $i=1$ and $i=2$,
then $\mu_\Gamma=0$ and $L^2(\Gamma\bk G)$ is tempered.
By Lemma~\ref{la:gifproduct}, $\psi_\Gamma\leq \rho$ in this case
which agrees with Corollary~\ref{cor:temperedness1}.
Moreover, for $\mu=(\mu_1,\mu_2)\in \overline{\mathfrak{a}^\ast_+}$
we have
\[
	\sup\|	\Re \wt\sigma_\Gamma \|_\mu = \max\left(\frac{\delta_{\Gamma_1}-\rho_1}{\mu_1},\frac{\delta_{\Gamma_2}-\rho_2}{\mu_2}, 0\right)
\]
and we see from Lemma~\ref{la:gifproduct} that this agrees with $\max(0,\delta_\Gamma'(\mu))$.
Let us next illustrate the implication of Theorem~\ref{thm:limit_cone_tempered_intro}:
Note that $\mathcal{L}_\Gamma = \overline{\mathfrak{a}_+}$
if $\Gamma_i$ are both infinite
and $\mathcal{L}_\Gamma\subseteq \mathfrak{a}_+\cup\{0\}$ if $\Gamma_i$ are both finite.
In the latter case $L^2(\Gamma\bk G)$ is tempered giving the same conclusion as Theorem~\ref{thm:limit_cone_tempered_intro}.
On the other hand, if $L^2(\Gamma \bk G)$ is non-tempered,
i.e.~if $\delta_{\Gamma_i} > \rho_i$ for at least one $i$,
then Theorem~\ref{thm:limit_cone_tempered_intro} implies that
$\mathcal{L}_\Gamma$ intersects the boundary of $\mathfrak{a}_+$
and the corresponding discrete subgroup must be infinite.

Finally, if $\Gamma_i$ is geometrically finite and non-cocompact in $\SL_2(\R)$
with $\delta_{\Gamma_i}-\rho_i >0$,
then $\delta_{\Gamma_i}-\rho_i$ is a discrete $L^2$-eigenvalue and
$[\rho_i^2,\infty[\subseteq \sigma(\Delta_i)$.
Hence, $(\delta_1-\rho_1,\delta_2-\rho_2)$ is a discrete joint
$L^2$ eigenvalue and there are also continuous spectral families on the boundaries:
In view of \eqref{eq:product_spectrum} this yields that there are continuous families of
joint spectra $(\pm(\delta_1-\rho_1), i\R) \in \widetilde\sigma_\Gamma$ and $(i\R,
\pm(\delta_2-\rho_2)) \in \widetilde\sigma_\Gamma$ which
lie on the boundary of the polyhedral region.

\subsection{\texorpdfstring{$\mathbf {\SL_3(\R)}$}{SL(3,R)} case}\label{ssec:sl3}
In the example $G=\SL_3(\R)$ or more generally if $G$ is locally isomorphic to 
$\mathfrak{sl}_3(\mathbb{K}), \mathbb{K}=\R,\C,\mathbb{H}$ or $\mathfrak{e}_{6(-26)}$,
the root system of restricted roots is $A_2$.
There are two simple roots $\alpha_1,\alpha_2$ with an angle of $2\pi/3$.
The half sum of positive roots $\rho$ is a multiple of the third positive root $\alpha_3=\alpha_1+ \alpha_2$.
For $G=\SL_3(\R)$ we have $\rho=\alpha_3$,
but more generally, if $m\coloneqq m_{\alpha_1}=m_{\alpha_2}=m_{\alpha_3}\geq 2$,
we have $\rho= m \alpha_3$.
For $\mathfrak{sl}_3(\C)$, $m=2$, for $\mathfrak{sl}_3(\mathbb{H})$, $m=4$, and for $\mathfrak{e}_{6(-26)}$, $m=8$.
The Weyl group consists of 6 elements, 3 rotations of an angle of $0, 2\pi/3$, $4\pi/3$, as well as 
the three reflections along the three positive roots.
Since $\widehat G_{\mathrm{sph}}\subseteq \{\lambda\in \mathfrak{a}_\C^\ast \mid -\overline{\lambda}\in W\lambda\}$,
for every $\lambda\in \widetilde \sigma_\Gamma$
with $\Re \lambda \neq 0$,
there is $i\in \{1,2,3\}$ with $\Re \lambda\in \R \alpha_i$ and $\Im \lambda \in \alpha_i^\perp$.

By $W$-invariance of $\widetilde \sigma_\Gamma$, 
we can always assume that $\Re \lambda\in \overline{\mathfrak{a}_+^\ast}$.
Hence, $\Re \lambda = r\rho$ with $r\geq 0$
and we note that $\|r\rho\|_{\mathrm{poly},\rho}=r$.
The general bound \eqref{eq:generalboundspectrum} implies $r\leq 1$
and the bound by Property (T) \eqref{eq:decay propT} implies $r\leq \frac{2m-1}{m}$
(which is $r\leq \frac 12$ for $G=\SL_3(\R)$).
\begin{figure}
	\centering
	\includegraphics[width=0.8\textwidth, trim= 2cm 2.3cm 3cm 2cm, clip]{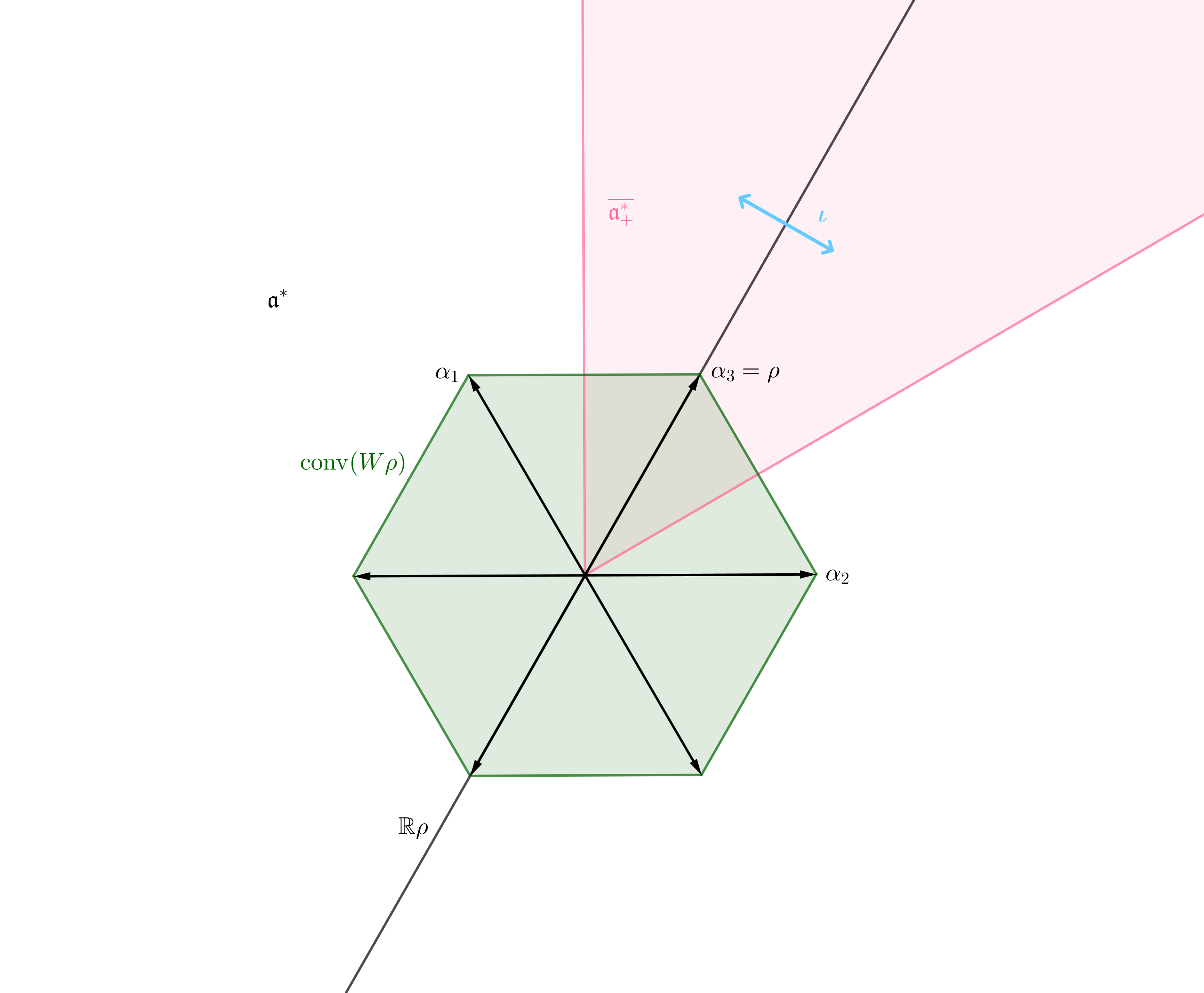}
	\caption{The root system of restricted roots for $\SL_3(\R)$. 
	The opposition involution $\iota$ is the reflection on the line spanned by $\rho$.}
	\label{fig:sl3basic}
\end{figure}

As $\mathfrak{a}^{\ast,\mathrm{Her}}\cap \overline{\mathfrak{a}^\ast_+} = \R_{\geq 0} \rho$,
we apply Theorem~\ref{thm:main} to $\rho$ and get
\[
\sup_{\lambda\in \widetilde \sigma_\Gamma} \|\Re \lambda\|_{\mathrm{poly},\rho} =
\max(0,\delta_{\Gamma}'(\rho)) \in [0,1].
\]
However, as mentioned in the introduction,
whenever $\delta_{\Gamma}'(\rho)>0$ we a priori neither know whether the supremum is attained,
nor, if it is attained, whether this supremum is real.
In this example however, we will be able to prove that both has to be the case.

Let us assume $\delta_{\Gamma}'(\rho)>0$ in the sequel. By definition and semi-continuity of $\psi_\Gamma$, there is $H_0\in \overline{\mathfrak{a}_+}$ such that $\psi_\Gamma (H_0)=(1+\delta_{\Gamma}'(\rho))\rho(H_0)$
and $\psi_\Gamma\leq (1+\delta_{\Gamma}'(\rho))\rho$.
However, $\psi_\Gamma$ and $\rho$ are invariant under the opposition involution
\[
 \iota(H)\coloneqq- w_0 H,\text{ where }w_0\in W\text{ is the longest Weyl group element}
\]
which is the negative of the reflection along $\alpha_3$, 
i.e.~the reflection on $\R\alpha_3 = \R\rho$.
Let us assume that $\Gamma$ is Zariski-dense, so that $\psi_\Gamma$ is concave (see Section~\ref{ss:gif}).
Therefore, $\frac 12 (H_0+\iota(H_0))\in \R\rho$ and $\psi_\Gamma(\frac 12(H_0+\iota(H_0))\geq \frac 12 \psi_\Gamma(H_0) + \frac 12 \psi_\Gamma(\iota(H_0))=\psi_\Gamma(H_0)$
by concavity and $\iota$-invariance of $\psi_\Gamma$.
Therefore, we can assume without loss of generality that $H_0\in \R\rho$.

Let us take a look at the bottom $\inf \sigma(\Delta)$ of the Laplace spectrum.
By \cite[Cor.~1.4]{WZ23} 
\[
\inf \sigma(\Delta)=|\rho|^2 - \max (0,\sup_{|H|=1} \psi_\Gamma(H) - \rho(H))^2.
\]
There exists $\lambda\in \widetilde \sigma_\Gamma$ with $\Re\lambda \in \mathfrak{a}_+^\ast\cup\{0\}$ and 
$\inf \sigma(\Delta)=\chi_\lambda (\Delta)=|\rho|^2- |\Re \lambda|^2 +|\Im\lambda|^2$.
Hence,
\begin{align*}
	|\Re \lambda|^2 
	&= |\Im \lambda|^2 + \max \left(0,\sup_{|H|=1} \psi_\Gamma(H) - \rho(H)\right)^2\\
	&\geq |\Im \lambda|^2 + \left(\frac{\delta_{\Gamma}'(\rho) \rho(H_0)}{|H_0|}\right)^2 \\
	&= |\Im \lambda|^2 + \delta_{\Gamma}'(\rho)^2 |\rho|^2.
\end{align*}
Here we used $\psi_\Gamma(H_0)-\rho(H_0)=\delta_{\Gamma}'(\rho) \rho(H_0)$ for the inequality
and $H_0\in \R\rho$ for the last equality.
On the other hand, since $\|\Re \lambda\|_{\mathrm{poly},\rho}\leq \delta_{\Gamma}'(\rho)$ we have
$|\Re \lambda|\leq \delta_{\Gamma}'(\rho)|\rho|$.
We conclude that
$\Im\lambda=0$ and $\Re\lambda = \delta_{\Gamma}'(\rho) \rho$,
i.e. $\delta_{\Gamma}'(\rho) \rho \in \widetilde \sigma_\Gamma$.

The case $\delta_{\Gamma}'(\rho) \leq 0$ means that $\psi_\Gamma\leq \rho$ and $\widetilde \sigma_\Gamma\subseteq i\mathfrak{a}^\ast$,
as well as $\inf \sigma(\Delta)=|\rho|^2$.
Therefore,
$\lambda\in \widetilde \sigma_\Gamma$ with $\chi_\lambda(\Delta)=|\rho|^2$ 
has to be $0=\delta_{\Gamma}'(\rho) \rho$.

To summarize, in the $A_2$ case with Zariski-dense $\Gamma$,
$\sup_{\lambda\in \widetilde \sigma_\Gamma} \|\Re \lambda\|_{\mathrm{poly}} = \theta$
is achieved at $\delta_{\Gamma}'(\rho) \rho$.
Let us emphasize, that our analysis provides no information
whether $\delta_{\Gamma}'(\rho) \rho$ is an isolated joint $L^2$-eigenvalue or is
part of continuous spectrum. However, as the joint spectral value $\delta_{\Gamma}'(\rho) \rho$
corresponds to the bottom  of the $L^2$-spectrum of $\Delta$ the recent work of
\cite{EFLO23} implies that (for Zariski dense $\Gamma$) $\delta_{\Gamma}'(\rho)\rho$ cannot
correspond to a joint $L^2$ eigenvalue of $\mathbb D(G/K)$ because otherwise
the bottom of the spectrum would be a $L^2$ eigenvalue contradicting their
result. We think that studying the properties of the spectrum inside the
polyhedral tubes $\delta_{\Gamma}'(\rho) \operatorname{conv}(W\rho) +i \mathfrak{a}^\ast$
is a highly interesting question which should be addressed in
the future, in particular as the $A_2$-case is the only higher rank case for which Theorem~\ref{thm:limit_cone_tempered_intro} does not apply.
  \small

  \subsection{Non-tempered example}\label{ssec:nontempered_exmpl}
  In this subsection, we illustrate how our polyhedral bounds yield precise knowledge on the location of exceptional spectrum for a new class of examples established by Fraczyk and Oh \cite{FO25}.

   They consider the following situation: Let $G=\mathrm{SO}_0(2,n)$, $n\geq
   3$, $H = \mathrm{SO}_0(1,n)\subset G$ and $\Gamma_0 < H$ a lattice.
   Then $\Gamma_0< G$ is a discrete subgroup for which the growth indicator
   function $\psi_{\Gamma_0}$ can be explicitly calculated and does not satisfy
   $\psi_{\Gamma_0}\leq\rho$. By construction, $\Gamma_0$ is not Zariski-dense.
   However, the bending techniques of \cite{Kas12} allow one to construct
   a family of discrete subgroups $(\Gamma_t)_{t>0}$, which, for
   sufficiently small $t$, are Zariski-dense and Anosov with respect to a
   maximal parabolic subgroup. Fraczyk and Oh then use that $\psi_{\Gamma_{t}}$
   depends in an appropriate sense continuously on $t$ and deduce
   non-temperedness of $\Gamma_t$ using Corollary~\ref{cor:temperedness1}. We
   now illustrate how the precise estimates on the polyhedral norms of
   Theorem~\ref{thm:main} not only show the existence of non-tempered
   spectrum, but also allow one to locate this part of the spectrum quite precisely. As a
   consequence, we will see that the examples in \cite{FO25} are optimal in the
   sense that the exceptional spectrum is arbitrarily close to the general
   bound enforced by Property (T) which is insurmountable for non-lattice
   subgroups.

  Let us introduce some notation in order to describe the spectrum:
  Let $G=\mathrm{SO}_0(2,n)$, $n\geq 3$, be the identity component of the indefinite special
  orthogonal group for a quadratic form of signature $(2,n)$.
  $\mathrm{SO}_0(2,n)$ has real rank $2$,
  so that $\mathfrak{a} = \R^2 \simeq \{\diag(v_1,v_2,0,\ldots,0,-v_2,-v_1)\colon v_1,v_2\in \R\}$.
  A choice of positive roots is
  \[
	  \Sigma^+ = \{\alpha_1\colon v\mapsto v_1-v_2,
		  \alpha_2\colon v\mapsto v_2,
		  \alpha_3 = \alpha_1+\alpha_2,
	  \alpha_4 = \alpha_1+2\alpha_2\},
  \]
  with multiplicities $m_{\alpha_1}=m_{\alpha_4}=1$ and $m_{\alpha_2}=m_{\alpha_3}=n-2$.
  Hence, $\mathfrak{a}_+ = \{v\mid v_1>v_2>0\}$
  and $\rho(v)=\frac 12 (nv_1 + (n-2)v_2)$.
  The root system is of type $B_2$, so that the Weyl group contains $-1$.
  In particular, the condition $-\overline{\lambda} \in W\lambda$ for
  $\lambda\in \widehat G_{\mathrm{sph}}$ gives no restriction for the real part.
  The quantitative bound on $\widehat G_{\mathrm{sph}}$ is given by
  $\Theta =\frac 12( \alpha_1+\alpha_4) \colon v\mapsto v_1$.
  More precisely, for any non-lattice discrete subgroup $\Gamma< \mathrm{SO}_0(2,n)$,
  \[
	  \psi_\Gamma(v)\leq (2\rho-\Theta)(v)=(n-1)v_1 + (n-2)v_2
  \]
  and
  \begin{equation}\label{eq:FO_prop_T_bound}
	  \Re \widetilde \sigma_\Gamma \subseteq \mathrm{conv}(W(\rho-\Theta)) = \frac {n-2}2 \mathrm{conv}(W\alpha_4),
  \end{equation}
  by \cite[Prop.~4.1]{FO25} and Theorem~\ref{thm:main}.
\begin{figure}
	\centering
	\includegraphics[width=\textwidth, trim= 0cm 0.5cm 0.5cm 0.5cm, clip]{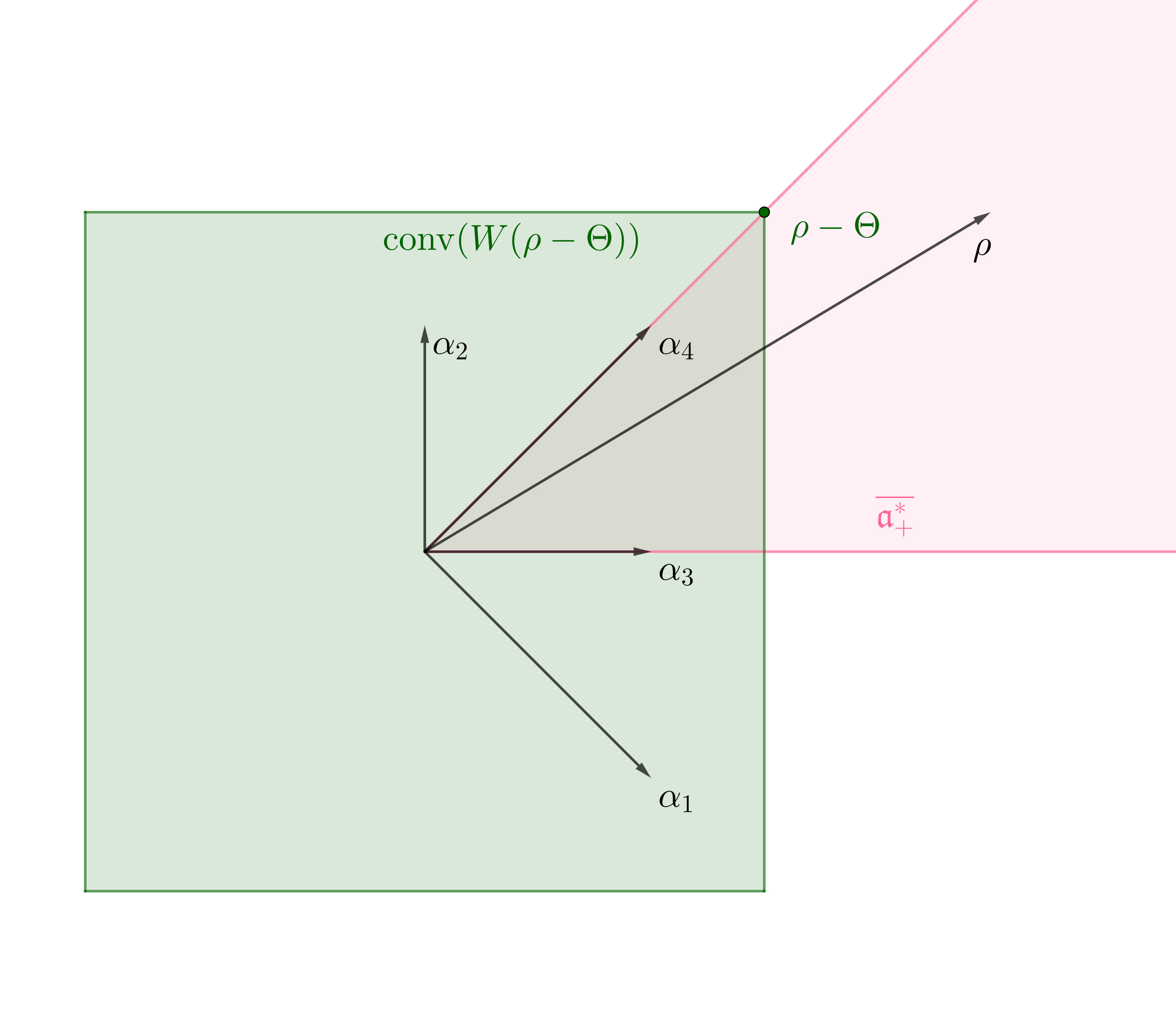}
	\caption{The root system of restricted roots for $\SO_0(2,5)$. 
	}
	\label{fig:B2}
\end{figure}

  Let $H=\mathrm{SO}_0(1,n)$ be identified with a subgroup of $\mathrm{SO}_0(2,n)$
  stabilizing a hyperplane on which the restriction of the quadratic form has signature $(1,n)$.
  We can choose the inclusion $H\to G$
such that $\mathfrak{h}\cap \mathfrak{a} = \ker \alpha_2$.
Let $\Gamma_0< \mathrm{SO}_0(1,n)$ be a lattice.
Then, by \cite[Prop.~4.2]{FO25}, we have
\[
	\psi_{\Gamma_0} (v) = \begin{cases}
		(n-1)v_1 &\colon v_1\geq 0, v_2=0,\\
		-\infty & \colon \text{else}
	\end{cases}
	.
\]
For $\mu(v)=\mu_1v_1 + \mu_2 v_2$ with $\mu_1\geq \mu_2\geq 0$,
$\mu\in \mathfrak{a}^{\ast,\mathrm{Her}}\cap \overline{\mathfrak{a}_+^\ast} = \overline{\mathfrak{a}_+^\ast}$
the polyhedral norm of $\lambda\in \overline{\mathfrak{a}^\ast_+}$
is $\|\lambda\|_{\mathrm{poly},\mu} = \sup_{v\in \mathfrak{a}_+} \frac{\lambda(v)}{\mu(v)}=
\sup_{1\geq v_2 \geq 0} \frac{\lambda_1 + \lambda_2v_2}{\mu_1+\mu_2v_2}$.
In the case $\mu_1=1$ and $\mu_2=0$, i.e.~$\mu=\alpha_3$,
$\|\lambda\|_{\mathrm{poly},\alpha_3} = \lambda_1+\lambda_2$,
and in the case $\mu_1=1$ and $\mu_2=1$, i.e.~$\mu=\alpha_4$,
\[
	\|\lambda\|_{\mathrm{poly},\alpha_4} = \sup_{1\geq v_2\geq 0} \frac{\lambda_1 +\lambda_2 v_2}{1+v_2} = 
	\sup_{v_2} \frac{\lambda_1-\lambda_2}{1+v_2} + \lambda_2=\lambda_1.
\]
The number $\delta_{\Gamma_0}'(\mu)$ equals
\[
	\delta_{\Gamma_0}'(\mu)=
	\sup_{v\in \overline{\mathfrak{a}_+}} \frac{\psi_{\Gamma_0}(v) - \rho(v)}{\mu(v)} =
	\frac{(n-1)-\frac n2}{\mu_1}= \frac{n-2}{2\mu_1}>0.
\]
and hence
\begin{equation}\label{eq:delta_Gamma_0_alpha_3}
 \delta_{\Gamma_0}'(\alpha_3) = \frac{n-2}{2}.
\end{equation}
Therefore,
\[
	\delta'_{\Gamma_0}(\mu)\mu(v) = \frac{n-2}2 \left(v_1 + v_2 \frac{\mu_2}{\mu_1}\right)
	\geq \frac{n-2}{2}v_1 = \frac{n-2}2 \alpha_3(v)
\]
for all $v\in \overline{\mathfrak{a}_+}$.
Hence,
\begin{equation}\label{eq:alpha_3_upper_bound}
\operatorname{conv}\left(W \frac{n-2}2 \alpha_3\right) \subseteq \delta_{\Gamma_0}'(\mu)\operatorname{conv}(W\mu)
\end{equation}
By Theorem~\ref{thm:main}, $\delta_{\Gamma_0}(\mu)= \sup\|\Re \wt \sigma_\Gamma\|_{\mathrm{poly},\mu}$ and with \eqref{eq:delta_Gamma_0_alpha_3} and \eqref{eq:alpha_3_upper_bound} imply $\mu_{\Gamma_0} = \frac{n-2}2 \alpha_3$.

Now \cite[Prop.~6.8]{FO25} yields that $\delta_{\Gamma_t}'(\mu)$ varies continuously in $t$\footnote{Actually,
\cite[Prop.~6.8]{FO25} only treats the case of $\mu=\rho$ but the proof applies verbatim for any $\mu\in \mathfrak{a}^\ast$
positive on $\overline{ \mathfrak{a}_+}$}.
Thus, for all $\varepsilon > 0$, there is $t>0$ such that $\delta_{\Gamma_{t}}(\mu)' \geq \frac{n-2}{2\mu_1}-\varepsilon$.
By Theorem~\ref{thm:main},
$\sup_{\lambda\in \widetilde \sigma_{\Gamma_{t}}} \|\Re \lambda\|_{\mathrm{poly},\mu} \geq
\frac{n-2}{2\mu_1} - \varepsilon$.
For $\mu=\alpha_4$ we obtain that there is a
$\lambda=\lambda(t, \varepsilon)\in \widetilde \sigma_{\Gamma_{t}}$
with $\|\Re \lambda\|_{\mathrm{poly},\alpha_4} \geq \frac{n-2}2-2\varepsilon$,
i.e.~$\Re \lambda_1 \geq \frac{n-2}{2}-2\varepsilon$.
On the other hand, for $\mu=\alpha_3$,
we have
$\Re \lambda_1 +\Re \lambda_2 = \|\Re \lambda\|_{\mathrm{poly},\alpha_3}\leq \delta_{\Gamma_{t}}'(\alpha_3) \leq \frac{n-2}2 +\varepsilon$.
Together this implies 
\[
	\left\|\Re \lambda- \frac{n-2}2\alpha_3\right\|= \mathcal{O}(\varepsilon).
\]
Similarly, one can show that $\|\mu_{\Gamma_t} - \frac{n-2}2 \alpha_3\| =\mathcal{O}(\varepsilon)$.

We have thus shown that the examples of Fraczyk and Oh allow one to construct discrete Zariski-dense subgroups of $\mathrm{SO}_0(2,n)$
with real part of the joint spectrum arbitrarily close to $\frac{n-2}{2}\alpha_3$.
As $\frac{n-2}{2}\alpha_3$ is exactly on the boundary of the region that bounds the exceptional spectrum for any non-lattice subgroup
by the quantitative Property (T) bound,
we see that their examples are in fact optimal,
in the sense that the exceptional spectrum for the Zariski-dense subgroups can become
as exceptional as possible under the general Property (T) bounds.
  
\begin{figure}
	\centering
	\includegraphics[width=\textwidth, trim= 0cm 1cm 0.5cm 0.5cm, clip]{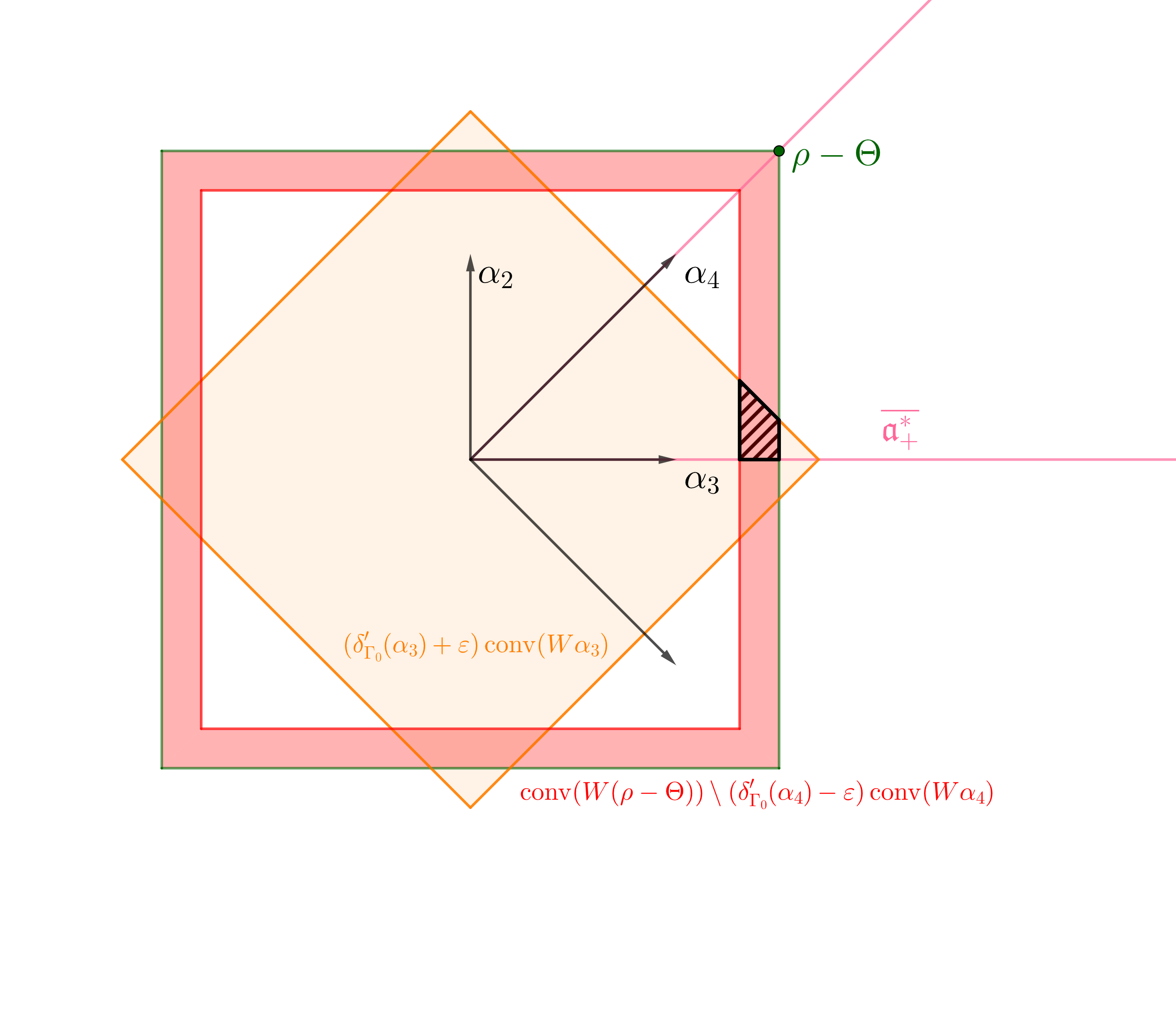}
	\caption{Using Theorem~\ref{thm:main} for $\alpha_4$ yields spectrum
		in the red region. 
		Applying it to $\alpha_3$ restricts to the orange region.
	Hence, by considering only $\overline{\mathfrak{a}_+^\ast}$,
	there must be spectrum in the black shaded region.
	}
	\label{fig:Nontemp}
\end{figure}

  \bibliographystyle{alpha}
  \bibliography{references_Revisions}

    \hrulefill

\end{document}